\documentclass[11pt]{article}
\usepackage{amssymb}
\usepackage{amsmath}
\usepackage{fancybox}
\usepackage{graphics}
\usepackage{latexsym}
\usepackage[applemac]{inputenc}
\usepackage{ae,aecompl}
\usepackage{amsthm}
\usepackage{mathrsfs}
\usepackage{amsfonts,amssymb}
 \usepackage[pdftex]{graphicx}

\DeclareGraphicsExtensions{.eps}

\newcommand{\op}[1]{\operatorname{#1 }}

\newcommand{\E}[1]{\mathbb{E} \left [ #1 \right ]}
\newcommand{\R}{\mathbb R}
\newcommand{\N}{\mathbb N}
\newcommand{\Z}{\mathbb Z}

\newcommand{\TT}{\mathbb T}

\renewcommand{\P}[2]{\mathbb{P}_{#2}\left( #1 \right)}

\def\be{{\bf e}}
\def\un{\underline}
\def\la{\longrightarrow}
\def\bm{{\bf m}}
\def\bM{{\bf M}}
\def\bp{{\bf p}}
\def\bq{{\bf q}}
\def\ov{\overline}
\def\eps{\varepsilon}

\def\build#1_#2^#3{\mathrel{
\mathop{\kern 0pt#1}\limits_{#2}^{#3}}}

\newtheorem{theorem}{Theorem}[section]

\newtheorem{proposition}[theorem]{Proposition}
\newtheorem{lemma}[theorem]{Lemma}
\newtheorem{corollary}[theorem]{Corollary}
\newtheorem{definition}[theorem]{Definition}

\newcommand{\m}{\mathcal{M}}
\renewcommand{\L}{\mathcal{L}}
\newcommand{\T}{\mathcal{T}}
\renewcommand{\P}[1]{\mathbb{P}\left[#1\right]}
\newtheorem{rek}[theorem]{Remark}

\usepackage{enumerate}
\def\llbracket{[\hspace{-.10em} [ }
\def\rrbracket{ ] \hspace{-.10em}]}
\title{\Huge{The Brownian cactus I\\ Scaling limits of discrete cactuses}}
\date{\tiny \today}
\author{Nicolas Curien, Jean-Fran\c cois Le Gall, Gr\'egory Miermont}
\begin{document}
\maketitle

\begin{abstract}
The cactus of a pointed graph is a discrete tree associated with this graph. Similarly,
with every pointed geodesic metric space $E$, one can associate an $\R$-tree called the
 continuous cactus of $E$.
We prove under general assumptions that the cactus of random planar maps
distributed according to Boltzmann weights and conditioned to have a fixed 
large number of vertices converges in distribution to a limiting space called the Brownian cactus,
in the Gromov-Hausdorff sense. Moreover, the Brownian cactus can be 
interpreted as the continuous cactus of the so-called Brownian map.
\end{abstract}

\section{Introduction}

In this work, we associate with every pointed graph a discrete tree 
called the cactus of the graph. Assuming that the pointed graph is chosen at random
in a certain class of planar maps with a given number of vertices,
and letting this number tend to infinity, we show that, modulo a suitable rescaling, the
associated cactus 
converges to a universal object, which we call the Brownian cactus.

In order to motivate our results, let us recall some basic facts
about planar maps.
A planar map is a proper embedding of 
a finite connected graph in the
two-dimensional sphere, viewed up to orientation-preserving
homeomorphisms of the sphere. The faces of the map are
the connected components of the complement of edges, and the degree of
a face counts the number of edges that are incident to it, with the convention that 
if both sides of an edge are incident to the same face, this edge is counted twice 
in the degree of the face.  Special cases
of planar maps are triangulations, where each face has degree $3$, quadrangulations,
where each face has degree $4$ and more generally $p$-angulations
where each face has degree $p$. Since the pioneering work of Tutte \cite{Tu}, planar maps have been 
thoroughly studied in
combinatorics, and they also arise in other areas of mathematics: See
in particular the book of Lando and Zvonkin \cite{LZ} for algebraic and 
geometric motivations.
Large
random planar graphs are of interest in theoretical
physics, where they serve as models of random geometry \cite{ADJ}.

A lot of recent work has been devoted to the study of scaling limits of
large random planar maps viewed as compact metric spaces. The vertex set
of the planar map is equipped with the graph distance, and one is interested in the
convergence of the (suitably rescaled) resulting metric space when the number
of vertices tends to infinity, in the sense of the Gromov-Hausdorff distance.
In the particular case of triangulations, this problem was stated by Schramm \cite{Sch}.
It is conjectured that, under mild 
conditions on the underlying distribution of the random planar map, this convergence holds
and the limit is the so-called Brownian map. Despite some recent
progress \cite{MM,MaMi,IM,BG0,Geo}, this conjecture is still
open, even in the simple case of uniformly distributed quadrangulations.
The main obstacle is the absence of a characterization
of the Brownian map as a random metric space. A compactness 
argument can be used to get the existence of sequential limits 
of rescaled random planar maps \cite{IM}, but the fact that there is
no available characterization of the limiting object prevents 
one from getting the desired convergence.

In the present work, we treat a similar problem, but we replace the
metric space associated with a planar map by a simpler metric space called
the cactus of the map. Thanks to this replacement, we are able to 
prove, in a very general setting,  the existence of a scaling limit,
which we call the Brownian cactus. Although this result remains far from
the above-mentioned conjecture, it gives another strong indication of the
universality of scaling limits of random planar maps, in the spirit of
the papers \cite{CS,MaMi,MieInvar,MW} which were concerned with the
profile of distances from a particular point.

Let us briefly explain the definition of the discrete cactus (see subsection \ref{sec:disccac}
for more details). We start from a graph $\mathbf{G}$ with a distinguished vertex $\rho$. 
Then, if $a$ and $b$ are two vertices of $\mathbf{G}$, and if $a_0=a,a_1,\ldots,a_p=b$
is a path from $a$ to $b$ in the graph $\mathbf{G}$, we consider the quantity
$$\op{d}_{\op{gr}}(\rho,a)+\op{d}_{\op{gr}}(\rho,b) - 2\min_{0\leq i\leq p} \op{d}_{\op{gr}}(\rho,a_i)$$
where $\op{d}_{\op{gr}}$ stands for the graph distance in $\mathbf{G}$. The cactus distance
$\op{d}_{\op{Cac}}^{\mathbf G}(a,b)$ is then the minimum of the preceding quantities
over all choices of a path from $a$ to $b$. The cactus distance is in fact only
a pseudo-distance: We have $\op{d}_{\op{Cac}}^{\mathbf G}(a,b)=0$ if and only if 
$\op{d}_{\op{gr}}(\rho,a)=\op{d}_{\op{gr}}(\rho,b) $ and if there is a path from $a$ to $b$ that stays at
distance at least $\op{d}_{\op{gr}}(\rho,a)$ from the point $\rho$. The cactus $\op{Cac}(\mathbf{G})$ associated with
$\mathbf{G}$ is the quotient space of the vertex set of $\mathbf{G}$ for the equivalence relation
$\asymp$ defined by putting $a\asymp b$ if and only if $\op{d}_{\op{Cac}}^{\mathbf G}(a,b)=0$. 
The set $\op{Cac}(\mathbf{G})$ is equipped by the distance induced by $\op{d}_{\op{Cac}}^{\mathbf G}$.
It is easy to verify that $\op{Cac}(\mathbf{G})$ is a discrete tree (Proposition \ref{cakd}).
Although much information is lost when going from $\mathbf{G}$ to its cactus, $\op{Cac}(\mathbf{G})$
still has a rich structure, as we will see in the case of planar maps.

A continuous analogue of the cactus can be defined for a (compact) geodesic metric space $\mathbf{E}$
having a distinguished point $\rho$. As in the discrete setting, the cactus distance between two points
$x$ and $y$ is the infimum over all continuous paths $\gamma$ from $x$ to $y$ of the difference
between the sum of the distances of $x$ and $y$ to the distinguished point $\rho$
and twice the minimal distance of a point of $\gamma$ to $\rho$. Again this is only a pseudo-distance, 
and the continuous cactus $\op{Kac}(\mathbf{E})$ is defined as the corresponding quotient space of $\mathbf{E}$. 
One can then check that the mapping $\mathbf{E}\la \op{Kac}(\mathbf{E})$ is continuous, and even 
Lipschitz, with respect to the Gromov-Hausdorff distance between pointed metric spaces 
(Proposition \ref{LipschitzGH}). It follows that if a sequence of (rescaled) pointed graphs $\mathbf{G}_n$ converges 
towards a pointed space $\mathbf{E}$ in the Gromov-Hausdorff sense, the (rescaled) cactuses $\op{Cac}(\mathbf{G}_n)$
also converge to $\op{Kac}(\mathbf{E})$. In particular, this implies that $\op{Kac}(\mathbf{E})$ is an $\R$-tree
(we refer to \cite{Evans} for the definition and basic properties of $\R$-trees).

The preceding observations yield a first approach to the convergence of rescaled cactuses
associated with random planar maps. Let $p\geq 2$ be an integer, and for every
$n\geq 2$, let $m_n$ be a random planar map that is uniformly distributed over the 
set of all rooted $2p$-angulations with $n$ faces
(recall that a planar map is rooted if there is a distinguished edge, which is oriented and whose origin
is called the root vertex). We view the vertex set $V(m_n)$ of $m_n$
as a metric space for the graph distance $\op{d}_{\op{gr}}$, with a distinguished point which is the
root vertex of the map. According to \cite{IM},
from any given strictly increasing sequence of
integers, we can extract a subsequence along which the rescaled pointed metric spaces 
$(V(m_n),n^{-1/4}\op{d}_{\op{gr}})$ converge in distribution in the Gromov-Hausdorff sense. As already 
explained above, the limiting distribution is not uniquely determined, and may depend
on the chosen subsequence. Still we call Brownian map any possible limit
that may arise in this convergence. Although the distribution of the Brownian map
has not been characterized, it turns out that the distribution of its continuous cactus 
is uniquely determined. Thanks to this observation, one easily gets that
the suitably rescaled discrete cactus of $m_n$ converges in distribution to a random
metric space (in fact a random $\R$-tree) which we call the Brownian cactus: See Corollary
\ref{conv-to-cactus} below.

Let us give a brief description of the Brownian cactus. The random $\R$-tree known as 
the CRT, which has been introduced and studied by Aldous \cite{Al1,Al3} is denoted by
$(\T_\be,\op{d}_\be)$. The notation $\T_\be$ refers to the fact that the CRT is conveniently
viewed as the $\R$-tree coded by a normalized Brownian excursion $\be=(\be_t)_{0\leq t\leq 1}$
(see Section \ref{Brcactus} for more details). 
Let $(Z_a)_{a\in\T_\be}$ be Brownian labels on the CRT. Informally, we may say that,
conditionally on $\T_\be$, $(Z_a)_{a\in\T_\be}$ is a centered Gaussian process 
which vanishes at the root of the CRT and satisfies $\mathbb{E}[(Z_a-Z_b)^2]=\op{d}_\be(a,b)$
for every $a,b\in\T_\be$. Let $a_*$ be the (almost surely unique) vertex of 
$\T_\be$ with minimal label. For every $a,b\in\T_\be$, let $\llbracket a,b\rrbracket$
stand for the geodesic segment between $a$ and $b$ in the tree $\T_\be$,
and set
$$\op{d}_{\op{KAC}}(a,b)= Z_a+Z_b - 2 \min_{c\in\llbracket a,b\rrbracket} Z_c.$$
Then $\op{d}_{\op{KAC}}$ is a pseudo-distance on $\T_\be$. The Brownian 
cactus $\op{KAC}$ is the quotient space of the CRT for this pseudo-distance. As explained above,
it can also be viewed as the continuous cactus associated with the Brownian map (here and later,
we abusively speak about ``the'' Brownian map although its distribution may not be unique).

The main result of the present work (Theorem \ref{convmapcactus}) states that the
Brownian cactus is also the limit in distribution of the discrete cactuses associated
with very general random planar maps. To explain this more precisely, we need 
to discuss Boltzmann distributions on planar maps. For technical reasons, we consider 
rooted and pointed planar maps, meaning that in addition to the root edge there is
a distinguished vertex. Let $\bq=(q_1,q_2,\ldots)$
be a sequence of non-negative weights satisfying general assumptions (we require 
that $\bq$ has finite support, that $q_k>0$ for some
$k\geq 3$, and that $\bq$ is critical in the sense of \cite{MaMi,MieInvar} -- the latter property can always be achieved by multiplying
$\bq$ by a suitable positive constant). For every rooted and pointed planar map $m$, set
$$W_\bq(m)= \prod_{f\in F(m)} q_{\op{deg}(f)}$$
where $F(m)$ stands for the set of all faces of $m$ and 
$\op{deg}(f)$ is the degree of the face $f$.
For every $n$, choose a random rooted and pointed planar map $M_n$ with 
$n$ vertices, in such a way that $\mathbb{P}(M_n=m)$
is proportional to $W_\bq(m)$ (to be precise, we need 
to restrict our attention to those integers $n$ such that 
there exists at least one planar map $m$ with $n$ vertices such that
$W_\bq(m)>0$). View $M_n$ as a graph pointed at
the distinguished vertex of $M_n$. Then Theorem \ref{convmapcactus}
gives the existence of a positive constant $B_\bq$ such that
$$B_\bq n^{-1/4}\cdot \op{Cac}(M_n)\build{\la}_{n\to\infty}^{(d)} \op{KAC}$$
in the Gromov-Hausdorff sense. Here the notation $\lambda \cdot E$
means that distances in the metric space $E$ are multiplied by
the factor $\lambda$.
This result applies in particular to uniformly distributed $p$-angulations
with a fixed number of faces (by Euler's formula the number of vertices
is then also fixed), and thus for instance to triangulations. In contrast with the
first approach described above, we do not need to restrict ourselves to
the bipartite case where 
$p$ is even.

As in much of the previous work on asymptotics for large random planar maps,
the proof of Theorem \ref{convmapcactus} relies on the existence \cite{BDG} of ``nice'' bijections
between planar maps and certain multitype labeled trees. It was observed
in \cite{MaMi} (for the bipartite case) and in \cite{MieInvar} that the tree 
associated with a random planar map following a Boltzmann distribution is
a (multitype) Galton-Watson tree, whose offspring distributions are determined 
explicitly in terms of the Boltzmann weights, and which is equipped with labels that
are uniformly distributed over admissible choices. This labeled tree can be
conveniently coded by the two random functions called the contour process
and the label process (see the end of subsection \ref{secshuffling}). In the bipartite case,
where $q_k=0$ if $k$ is odd, one can prove \cite{MaMi} that the contour process and the label process
associated with the random planar map $M_n$ converge as $n\to\infty$, modulo
a suitable rescaling, towards the pair consisting of a normalized Brownian excursion
and the (tip of the) Brownian snake driven by this excursion. This convergence is
a key tool for studying  the convergence of rescaled (bipartite) random planar maps
towards the Brownian map \cite{IM}. In our general non-bipartite setting, it 
is not known whether the preceding convergence still holds, but Miermont \cite{MieInvar} observed
that it does hold if the tree is replaced by a ``shuffled'' version.
Fortunately for our purposes, although the convergence of the coding functions 
of the shuffled tree would not be effective
to study the asymptotics of rescaled planar maps, it gives enough
information to deal with the associated cactuses. This is one of the
key points of the proof of Theorem \ref{convmapcactus} in Section 
\ref{convcactusmaps}.

The last two sections of the present work are devoted to some properties of
the Brownian cactus. We first show that the Hausdorff dimension of
the Brownian cactus is equal to $4$ almost surely, and is therefore the same
as that of the Brownian map computed in \cite{IM}. As a tool for the calculation of
the Hausdorff dimension, we derive precise information on the volume of
balls centered at a typical point of the Brownian cactus (Proposition \ref{haus}).
Finally, we apply  ideas of the theory of the Brownian cactus to a
problem about the geometry of the Brownian map. Precisely, given 
three ``typical'' points in the Brownian map, we study the existence and uniqueness of a 
cycle with minimal length that separates the first point from the second one and visits the third one.
This is indeed a continuous version of a problem discussed by Bouttier and 
Guitter \cite{BG} in the 
discrete setting of large quadrangulations. In particular, we recover the explicit 
distribution of the volume of the connected components bounded by the minimizing
cycle, which had been derived in \cite{BG} via completely different methods.
The results of this section strongly rely on the study of geodesics in the
Brownian map developed in \cite{Geo}.

The subsequent paper \cite{cactus2} derives further results about the Brownian
cactus and in particular studies the asymptotic behavior of the number of 
``branches'' of the cactus above level $h$ that hit level $h+\eps$, when
$\eps$ goes to $0$. In terms of the Brownian map, if 
$B(\rho,h)$ denotes the open ball of radius $h$ centered at the root $\rho$ and 
$N_{h,\eps}$ denotes the
number of connected components of the complement 
of $B(\rho,h)$ that intersect the complement of $B(\rho,h+\eps)$, the main
result of \cite{cactus2} states that $\eps^3N_{h,\eps}$ converges as 
$\eps$ goes to $0$ to  a nondegenerate random variable. This convergence is
closely related to an upcrossing approximation for the local time of
super-Brownian motion, which is of independent interest. 

The paper is organized as follows. In Section \ref{cactusDC}, we give the definitions
and main properties
of discrete and continuous cactuses, and establish connections between the discrete
and the continuous case. In Section \ref{Brcactus}, after recalling the construction and
main properties of the
Brownian map, we introduce the Brownian cactus and show that it coincides 
with the continuous cactus of the Brownian map. Section \ref{convmapcactus} 
contains the statement and the proof of our main result Theorem \ref{convmapcactus}.
As a preparation for the proof, we recall in subsection \ref{biject} the construction and
main properties of the bijections between planar maps and 
multitype labeled trees. Section \ref{hauscalc} is devoted to the Hausdorff dimension
of the Brownian cactus, and Section \ref{cycle} deals with minimizing cycles in
the Brownian map. An appendix gathers some facts about planar maps
with Boltzmann distributions, that are needed in Section
\ref{convcactusmaps}. 

\smallskip
\noindent{\bf Acknowledgement.} We thank Itai Benjamini for the name cactus as well
as  for suggesting the study of this mathematical object.

\section{Discrete and continuous cactuses} 
\label{cactusDC}

\subsection{The discrete cactus}
\label{sec:disccac}

 Throughout this section, we consider a graph ${G}=(V,\mathcal{E})$, meaning that $V$ is a finite set
 called the vertex set and $\mathcal{E}$ is a subset of the set of all (unordered) pairs $\{v,v'\}$ of distinct 
 elements of $V$.

If $v,v' \in V$, a \emph{path} from $v$ to $v'$  in ${G}$ is a finite sequence $\gamma = (v_{0}, \ldots , v_{n})$ in $V$,
such that $v_{0}=v$, $v_{n}=v'$ and $\{v_{i-1},v_{i}\}\in \mathcal{E}$, for every $1\leq i \leq n$. The integer $n\geq 0$ is called the length
of $\gamma$. We assume that ${G}$ is connected, so that a path from $v$ to $v'$ exists for every choice of 
$v$ and $v'$. 
The \emph{graph distance} $\op{d}^{{G}}_{\op{gr}}(v,v')$ is the minimal length of a path from $v$ to $v'$ in $G$. A path with minimal
length is called a geodesic from $v$ to $v'$ in $G$.

In order to define the cactus distance we consider also
a distinguished 
 point $\rho$  in $V$. The triplet 
$\mathbf{G}=(V,\mathcal{E},\rho)$ is then called a pointed graph. With this pointed graph we associate
the cactus (pseudo-)distance defined by setting for  every $v,v' \in V$,
$$ \op{d_{Cac}^{\mathbf{G}}}(v,v') := \op{d}_{\op{gr}}^{{G}}(\rho,v) + \op{d}_{\op{gr}}^{{G}}(\rho,v') -2 \max_{\gamma : v \to v'}\min_{a \in \gamma} \op{d}_{\op{gr}}^{{G}}(\rho,a),$$
where the maximum is over all paths $\gamma$ from $v$ to $v'$ in ${{G}}$. 

\begin{proposition} \label{cactusdistance} The mapping  $(v,v')\to\op{d_{Cac}^\mathbf{G}}(v,v')$ is a pseudo-distance on $V$ taking integer values. Moreover, for every $v,v'\in V$,
\begin{equation}
\label{uppercactus-graph}
\op{d}_{\op{gr}}^{G}(v,v') \geq \op{d_{Cac}^\mathbf{G}}(v,v').\end{equation}
and 
\begin{equation}
\label{dist-root}
\op{d_{Cac}^\mathbf{G}}(\rho,v)=\op{d}_{\op{gr}}^{G}(\rho,v).
\end{equation}
\end{proposition}

\proof It is obvious that $\op{d_{Cac}^\mathbf{G}}(v,v)=0$
and $\op{d_{Cac}^\mathbf{G}}(v,v')=\op{d_{Cac}^\mathbf{G}}(v',v)$. Let us verify the triangle inequality. Let $v,v',v'' \in V$ and choose two paths $\gamma_{1}:v\to v'$ and $\gamma_{2}:v'\to v''$  such that $\min_{a\in \gamma_1} \op{d}_{\op{gr}}^{G}(\rho,a)$ is maximal among all paths $\gamma : v \to v'$ in ${G} $ and a similar property holds for $\gamma_{2}$. The concatenation of $\gamma_{1}$ and $\gamma_{2}$ gives a path $\gamma_{3}:v \to v''$ and we easily get 
$$ \op{d}_{\op{Cac}}^\mathbf{G}(v,v'') \leq \op{d}_{\op{gr}}^{G}(\rho,v) + \op{d}_{\op{gr}}^{G}(\rho,v'') -2 \min_{a \in \gamma_{3}}\op{d}_{\op{gr}}^ {G}(\rho,a) \leq \op{d}_{\op{Cac}}^\mathbf{G}(v,v') + \op{d}_{\op{Cac}}^\mathbf{G}(v',v'').$$
In order to get the bound $\eqref{uppercactus-graph}$, let $v,v'\in V$, and choose a geodesic path $\gamma$ from $v$ to $v'$. Let $w$ be a point on the path $\gamma$ whose distance to $\rho$ is minimal. Then,
\begin{align*}
\op{d}^{G}_{\op{gr}}(v,v') = \op{d}^{G}_{\op{gr}}(v,w) + \op{d}^{G}_{\op{gr}}(w,v')
&\geq \op{d}^{G}_{\op{gr}}(\rho,v)+ \op{d}^{G}_{\op{gr}}(\rho,v') - 2\op{d}^{G}_{\op{gr}}(\rho,w)\\
&=\op{d}^{G}_{\op{gr}}(\rho,v)+ \op{d}^{G}_{\op{gr}}(\rho,v') -  2\min_{a\in \gamma} \op{d}^{G}_{\op{gr}}(\rho,a)\\
&\geq \op{d}_{\op{Cac}}^\mathbf{G}(v,v').
\end{align*}
Property (\ref{dist-root}) is immediate from the definition.
\endproof

As usual, we introduce the equivalence relation $\overset{\mathbf{G}}{\asymp}$ defined on $V$ by setting 
$v \overset{\mathbf{G}}{\asymp} v'$ if and only $ \op{d_{Cac}^\mathbf{G}}(v,v')=0$. Note that $v \overset{\mathbf{G}}{\asymp} v'$
if and only if $\op{d}_{\op{gr}}^{G}(\rho,v)=\op{d}_{\op{gr}}^{G}(\rho,v')$ and there exists a path from $v$ to $v'$ 
that stays at distance at least $\op{d}_{\op{gr}}^{G}(\rho,v)$ from $\rho$.

The corresponding 
quotient space is denoted by   $\op{Cac}(\mathbf{G})= V\,/ \overset{\mathbf{G}}{\asymp}$. The pseudo-distance $\op{d_{Cac}^\mathbf{G}}$
induces a distance on $\op{Cac}(\mathbf{G})$, and we keep the notation 
$\op{d_{Cac}^\mathbf{G}}$ for this distance.

\begin{proposition}\label{cakd} Consider the graph $G^\circ$ whose vertex set is $V^\circ=\op{Cac}(\mathbf{G})$ and whose edges are all pairs
$\{a,b\}$ such that $\op{d}_{\op{Cac}}^\mathbf{G}(a,b)=1$. Then this graph is a tree, and the graph distance 
$\op{d}^{G^\circ}_{\op{gr}}$ on $V^\circ$ coincides with the cactus distance $\op{d_{Cac}^\mathbf{G}}$ on $\op{Cac}(\mathbf{G})$.
\end{proposition}

\proof Let us first verify that the graph $G^\circ$ is a tree. If $u \in V$ we use the notation $\overline{u}$ for the equivalence class of $u$
in the quotient  $\op{Cac}(\mathbf{G})$. We argue by contradiction and assume that there exists a 
(non-trivial) cycle in $\op{Cac}(\mathbf{G})$. We can then find an integer $n\geq 3$ and
vertices  $x_{0},x_{1},x_{2}, \ldots ,x_{n} \in V$ such that 
$$ \left \{ \begin{array}{l} \overline{x}_{0}= \overline{x}_{n}, \\ \op{d}_{\op{Cac}}^\mathbf{G}(x_{i},x_{i+1})=1, \mbox{ for every } 0\leq i \leq n-1,\\
 \overline{x}_{0}, \overline{x}_{1}, \ldots ,\overline{x}_{n-1} \mbox{ are distinct}. \end{array}\right.$$  

Without loss of generality, we may assume that $\op{d}_{\op{gr}}^{G}(\rho,x_{0})=\max\{ \op{d}_{\op{gr}}^{G}(\rho,x_{i}), 0 \leq i \leq n\}.$ By \eqref{dist-root}, we have $|\op{d}_{\op{gr}}^{G}(\rho,x_{0})-\op{d}_{\op{gr}}^{G}(\rho,x_{1})|\leq 
\op{d}_{\op{Cac}}^\mathbf{G}(x_{0},x_{1})=1$.
If $\op{d}_{\op{gr}}^{G}(\rho,x_{0})= \op{d}_{\op{gr}}^{G}(\rho,x_{1})$ then it follows from the definition of $\op{d}_{\op{Cac}}^\mathbf{G}$ that $\op{d}_{\op{Cac}}^\mathbf{G}(x_{0},x_{1})$ is even and thus different from $1$. So we must have $$\op{d}_{\op{gr}}^{G}(\rho,x_{1})=\op{d}_{\op{gr}}^{G}(\rho,x_{0})-1.$$ Combining this equality with
the property $\op{d}_{\op{Cac}}^\mathbf{G}(x_{0},x_{1})=1$, we obtain that there exists a path from $x_{0}$ to $x_{1}$ that stays at distance at least $\op{d}_{\op{gr}}^{G}(\rho,x_1)$ from $\rho$.

Using the same arguments and the equality $\op{d}_{\op{Cac}}^\mathbf{G}(x_{0},x_{n-1})=1$, we obtain similarly
that $\op{d}_{\op{gr}}^{G}(\rho,x_{n-1})=\op{d}_{\op{gr}}^{G}(\rho,x_{0})-1=\op{d}_{\op{gr}}^{G}(\rho,x_1)$
and that there exists a path from $x_{n-1}$ to $x_0$ that stays at distance at least $\op{d}_{\op{gr}}^{G}(\rho,x_{1})$ from $\rho$.

Considering the  concatenation of the two paths we have constructed, we get $\op{d}_{\op{Cac}}^\mathbf{G}(x_{1},x_{n-1})=0$ or equivalently $\overline{x}_{1}=\overline{x}_{n-1}$. This gives the desired contradiction, and we have proved that $G^\circ$ is a tree. 

We still have to verify the equality of the distances $\op{d}^{G^\circ}_{\op{gr}}$ and $\op{d}_{\op{Cac}}^\mathbf{G}$ on
$\op{Cac}(\mathbf{G})$. The bound $\op{d}_{\op{Cac}}^\mathbf{G}\leq \op{d}^{G^\circ}_{\op{gr}}$ is immediate from
the triangle inequality for $\op{d}_{\op{Cac}}^\mathbf{G}$ and the existence of a geodesic between any 
pair of vertices of $G^\circ$. Conversely, let $a,b\in\op{Cac}(\mathbf{G})$. We can find a path $(y_0,y_1,\ldots,y_n)$
in $G$ such that $\overline{y}_0=a$, $\overline{y}_n=b$ and
$$\op{d}_{\op{Cac}}^\mathbf{G}(a,b)= \op{d}^{G}_{\op{gr}}(\rho,y_0) +  \op{d}^{G}_{\op{gr}}(\rho,y_n) 
-2 \min_{0\leq j\leq n}  \op{d}^{G}_{\op{gr}}(\rho,y_j).$$
Put $m=\min_{0\leq j\leq n}  \op{d}^{G}_{\op{gr}}(\rho,y_j)$, $p=\op{d}^{G}_{\op{gr}}(\rho,y_0)$
and $q=\op{d}^{G}_{\op{gr}}(\rho,y_n) $ to simplify notation. Then set, for every $0\leq i\leq p-m$,
$$k_i=\min\{j\in\{0,1,\ldots,n\}: \op{d}^{G}_{\op{gr}}(\rho,y_j)= p-i\}$$
and, for every $0\leq i\leq q-m$,
$$\ell_i=\max\{j\in\{0,1,\ldots,n\}: \op{d}^{G}_{\op{gr}}(\rho,y_j)=q-i\}.$$
Then $\ov y_{k_0},\ov y_{k_1},\ldots,\ov y_{k_{p-m}}=\ov y_{\ell_{q-m}},\ov y_{\ell_{q-m-1}},\ldots, \ov y_{\ell_1},\ov y_{\ell_0}$ is a path
from $a$ to $b$ in
$G^\circ$. It follows that
$$\op{d}^{G^\circ}_{\op{gr}}(a,b)\leq p+q-2m = \op{d}_{\op{Cac}}^\mathbf{G}(a,b),$$
which completes the proof.
\endproof

\begin{rek} {\rm The notion of the cactus associated with a pointed graph strongly depends on the
choice of the distinguished point $\rho$.}
\end{rek}

In the next sections, we will be interested in rooted planar maps, which will even be pointed in
Section \ref{convcactusmaps}. With such a planar map, we can associate a pointed graph in the
preceding sense: just say that $V$ is the vertex set of the map, $\mathcal{E}$ is the set of all
pairs $\{v,v'\}$ of distinct points of $V$ such that there exists (at least) one edge of the map 
between $v$ and $v'$, and the vertex $\rho$ is either the root vertex, for a map that is
only rooted, or the distinguished point for a map that is rooted and pointed. Note that the
graph distance corresponding to this pointed graph (obviously) coincides with the 
usual graph distance on the vertex set of the map. Later, when we speak about the
cactus of a planar map, we will always refer to the cactus of the associated pointed graph.
In agreement with the notation of this section, we will use bold letters $\bm,\bM$
to denote the pointed graphs associated with the planar maps $m,M$. 

\begin{figure}[h]
\begin{center}
\includegraphics[width=16cm]{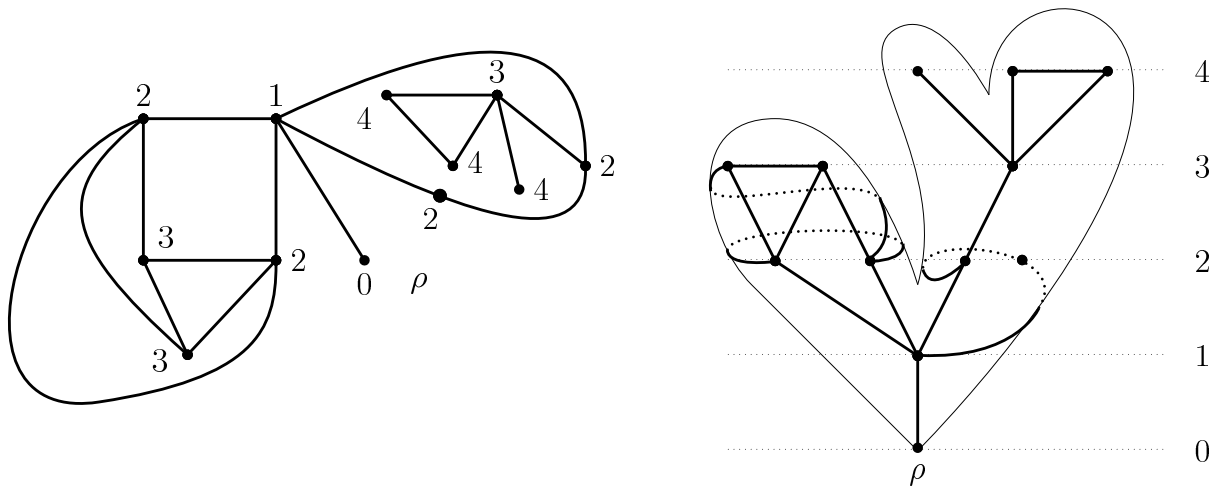}
\caption{A planar map and on the right side the same planar map represented so that the height of every vertex coincides with its
distance from the distinguished vertex $\rho$. We see a tree structure emerging from this picture, which corresponds to the associated cactus. }
\end{center}
\end{figure}

 \subsection{The continuous cactus} 
 \label{sec:contiK}
 
Let us recall some basic notions from metric geometry. If $(E,d)$ is a metric space and 
$\gamma:[0,T] \longrightarrow E$ is a continuous curve in $E$, the length of $\gamma$
is defined by:
$$L(\gamma)=\sup_{0=t_{0}< \cdots<t_{k}=T}\sum_{i=0}^{k-1}d\big(\gamma(t_{i}),\gamma(t_{i+1})\big),$$
where the supremum is over all choices of the subdivision $0=t_0<t_1<\cdots<t_k=T$ of $[0,T]$. 
Obviously $L(\gamma)\geq d(\gamma(0),\gamma(T))$. 

We say that $(E,d)$ is a geodesic space if for every $a,b\in E$ there exists a continuous curve
$\gamma:[0,d(a,b)]\longrightarrow E$ such that $\gamma(0)=a$, $\gamma(d(a,b))=b$ and $d(\gamma(s),\gamma(t))=t-s$
for every $0\leq s\leq t\leq d(a,b)$. 
Such a curve $\gamma$ is then called a geodesic from $a$ to $b$ in $E$. Obviously, $L(\gamma)=d(a,b)$. 
A pointed geodesic metric space is a geodesic space with a distinguished point $\rho$. 
 
Let ${\mathbf E}=(E,d,\rho)$ be a pointed geodesic compact metric space. We define the (continuous) cactus 
associated with $(E,d,\rho)$ in a way very similar to what we did in the discrete setting. We first
define for every $a,b\in E$,
$$\op{d}^{\mathbf E}_{\op{Kac}}(a,b)= d(\rho,a)+d(\rho,b) - 2 \sup_{\gamma:a\to b} \Big(\min_{0\leq t\leq 1} d(\rho,\gamma(t))\Big),$$
where the supremum is over all continuous curves $\gamma:[0,1]\longrightarrow E$
such that $\gamma(0)=a$ and $\gamma(1)=b$.

The next proposition is then analogous to Proposition \ref{cactusdistance}.

\begin{proposition}
\label{contcactus}
The mapping $(a,b)\longrightarrow \op{d}^{\mathbf E}_{\op{Kac}}(a,b)$ is a pseudo-distance on $E$. Furthermore, for
every $a,b\in E$,
$$\op{d}^{\mathbf E}_{\op{Kac}}(a,b)\leq d(a,b)$$
and 
$$\op{d}^{\mathbf E}_{\op{Kac}}(\rho,a)=d(\rho,a).$$
\end{proposition}

The proof is exactly similar to that of Proposition \ref{cactusdistance}, and we leave the details to the 
reader. Note that in the proof of the bound $\op{d}^{\mathbf E}_{\op{Kac}}(a,b)\leq d(a,b)$ we use the existence of a geodesic from
$a$ to $b$. 

If $a,b\in E$, we put $a\overset{\mathbf E}{\asymp} b$ if
$\op{d}^{\mathbf E}_{\op{Kac}}(a,b)=0$.  We define the cactus of
$(E,d,\rho)$ as the quotient space $\op{Kac}(\mathbf{E}):=
E\,/\overset{\mathbf E}{\asymp} $, which is equipped with the quotient
distance $\op{d}^{\mathbf E}_{\op{Kac}}$. Then $\op{Kac}(\mathbf{E})$
is a compact metric space, which is pointed at the equivalence class
of $\rho$. 

\begin{rek} 
\label{minimiz} {\rm It is natural to ask whether the supremum in the definition of $\op{d}^{\mathbf E}_{\op{Kac}}(a,b)$
is achieved, or equivalently whether there is a continuous path $\gamma$ from $a$ to $b$ such that
$$\op{d}^{\mathbf E}_{\op{Kac}}(a,b)= d(\rho,a)+d(\rho,b) - 2\min_{0\leq t\leq 1} d(\rho,\gamma(t)).$$
We will return to this question later.}
\end{rek}

\subsection{Continuity properties of the cactus}
 Let us start by recalling the definition of  the Gromov-Hausdorff distance between two
 pointed compact metric spaces (see \cite{Gro} and \cite[Section 7.4]{BBI01} for more details).
 
Recall that if $A$ and $B$ are two compact subsets of a metric space $(E,d)$, the Hausdorff distance between $A$ and $B$ is 
$$ \op{d}_{\op{H}}^E(A,B) := \inf\{ \varepsilon>0 : A \subset B^\varepsilon \mbox{ and }B \subset A^\varepsilon\},$$
where $X^\varepsilon := \{ x \in E : d(x,X) \leq \varepsilon\}$ denotes the $\varepsilon$-neighborhood of a subset $X $ of $E$.

\begin{definition}
If $\mathbf{E}=(E,d,\rho)$ and $\mathbf{E'}=(E,d',\rho')$ are two pointed compact  metric spaces, the Gromov-Hausdorff distance between
$\mathbf{E}$ and $\mathbf{E'}$  is 
$$ \op{d_{GH}}(\mathbf{E},\mathbf{E'}) = \inf\big\{\op{d}_{\op{H}}^F(\phi(E),\phi'(E')) \vee \delta(\phi(\rho),\phi'(\rho')) \big\},$$ where the infimum is taken over all choices of the metric space $(F,\delta)$ and the  isometric embeddings $\phi : E \to F$ and $\phi'  : E' \to F$ of $E$ and $E'$ into $F$.  \end{definition}
 The  Gromov-Hausdorff  distance is indeed a metric on the space of isometry classes of pointed compact metric spaces. An alternative definition of this distance uses \textit{correspondences.} A correspondence between two pointed metric spaces $(E,d,\rho)$ and $(E',d',\rho')$ is a subset $\mathcal{R}$ of $E\times E'$ containing $(\rho,\rho')$,  such that, for every $x_{1} \in E$, there exists at least one point $x_{2}\in E'$ such that $(x_{1},x_{2}) \in \mathcal{R}$ and conversely, for every $y_{2}\in E'$, there exists at least one point $y_{1}\in E$ such that $(y_{1},y_{2}) \in \mathcal{R}$. The distortion of the correspondence $\mathcal{R}$ is defined by 
 $$ \op{dis}(\mathcal{R}) := \sup_{}\big\{|d(x_{1},y_{1})-d'(x_{2},y_{2})| : (x_{1},x_{2}),(y_{1},y_{2}) \in \mathcal{R} \big\}.$$
The Gromov-Hausdorff distance can be expressed in terms of correspondences by the formula
 \begin{equation}
 \label{GHcorres}
 \op{d_{GH}}(\mathbf{E},\mathbf{E'})= \frac{1}{2} \inf \big\{\hspace{-0.5mm}\op{dis}(\mathcal{R})\big\},
 \end{equation} where the infimum is over all correspondences $\mathcal{R}$ between $\mathbf{E}$ and $\mathbf{E'}$. 
 See \cite[Theorem 7.3.25]{BBI01} for a proof in the non-pointed case, which is easily adapted.

\begin{proposition}
\label{LipschitzGH}
Let ${\mathbf E}$ and $\mathbf{E'}$ be two pointed geodesic compact metric spaces. Then,
$$ \op{d_{GH}}(\op{Kac}(\mathbf{E}),\op{Kac}(\mathbf{E'}))\leq 6\,\op{d_{GH}}(\mathbf{E},\mathbf{E'}).$$
\end{proposition}

\proof 
It is enough to
verify that, for any correspondence ${\mathcal R}$ between $\mathbf{E}$ and $\mathbf{E'}$ with distortion
$D$, we can find a correspondence ${\mathscr R}$ between $\op{Kac}(\mathbf{E})$ and $\op{Kac}({\mathbf{E'}})$ 
whose distortion is bounded above by $6D$. We define ${\mathscr R}$ as the set of all
pairs $(a,a')$ such that there exists (at least) one representative $x$ of $a$ in $E$ and one representative
$x'$ of $a'$ in $E'$, such that $(x,x')\in\mathcal{R}$. 

Let $(x,x')\in \mathcal{R}$ and $(y,y')\in \mathcal{R}$. We need to verify that
$$|\op{d}^{\mathbf E}_{\op{Kac}}(x,y)- \op{d}^{\mathbf E'}_{\op{Kac}}(x',y')| \leq 6D.$$
Fix $\varepsilon>0$. We can find a continuous curve $\gamma:[0,1]\longrightarrow E$ such that
$\gamma(0)=x$, $\gamma(1)=y$ and 
$$
d(\rho,x)+d(\rho,y) -2\min_{0\leq t\leq 1} d(\rho,\gamma(t)) \leq \op{d}^{\mathbf E}_{\op{Kac}}(x,y) + \varepsilon.$$
By continuity, we may find a subdivision $0=t_0<t_1<\cdots<t_p=1$
of $[0,1]$ such that $d(\gamma(t_i),\gamma(t_{i+1}))\leq D$ for every
$0\leq i\leq p-1$. For every $0\leq i\leq p$, put $x_i=\gamma(t_i)$, and choose $x'_i\in E'$ such that
$(x_i,x'_i)\in\mathcal{R}$. We may and will take $x'_0=x'$ and $y'_0=y'$. Now note that, for
$0\leq i\leq p-1$,
$$d'(x'_i,x'_{i+1})\leq d(x_i,x_{i+1}) + D \leq 2D.$$
Since $\mathbf{E'}$ is a geodesic space, we can find a curve $\gamma':[0,1]\longrightarrow E'$
such that $\gamma'(t_i)=x'_i$, for every $0\leq i\leq p$, and any point $\gamma'(t)$, $0\leq t\leq 1$
lies within distance at most $D$ from one of the points $\gamma'(t_i)$. It follows that
$$\min_{0\leq t\leq 1}d'(\rho',\gamma'(t))\geq \min_{0\leq i\leq p}d'(\rho',\gamma'(t_i))-D
\geq \min_{0\leq i\leq p}d(\rho,\gamma(t_i))-2D.$$
Hence,
\begin{align*}
\op{d}^{\mathbf E'}_{\op{Kac}}(x',y')&\leq d'(\rho',x')+d'(\rho',y') - 2 \min_{0\leq t\leq 1}d'(\rho',\gamma'(t))\\
&\leq d(\rho,x) + d(\rho,y) - 2\min_{0\leq t\leq 1}d(\rho,\gamma(t))+ 6D\\
&\leq \op{d}^{\mathbf E}_{\op{Kac}}(x,y) + 6D + \varepsilon
\end{align*}
The desired result follows since
$\varepsilon$ was arbitrary and  we can interchange the roles of $\mathbf{E}$ and $\mathbf{E'}$. 
\endproof

\subsection{Convergence of discrete cactuses}

Let $\mathbf{G}=(V,\mathcal{E},\rho)$ be a pointed graph (and write $G=(V,
\mathcal{E})$ for the non-pointed graph as previously). We can identify $\mathbf{G}$ with the
pointed (finite) metric space $(V,\op{d}^{G}_{\op{gr}},\rho)$. For any 
real $r>0$, we then denote the ``rescaled graph''  $(V,r \op{d}^{G}_{\op{gr}},\rho)$
by $r\cdot \mathbf{G}$.

Similarly, we defined $\op{Cac}(\mathbf{G})$ as a pointed finite metric space. The
space $r\cdot \op{Cac}(\mathbf{G})$ is then obtained by multiplying the
distance on $\op{Cac}(\mathbf{G})$ by the factor $r$.

\begin{proposition}
\label{convdiscactus}
Let $(\mathbf{G}_n)_{n\geq 0}$ be a sequence of pointed graphs, and let $(r_n)_{n\geq 0}$
be a sequence of positive real numbers converging to $0$. Suppose that 
$r_n\cdot \mathbf{G}_n$ converges to a pointed compact metric space $\mathbf{E}$, in the
sense of the Gromov-Hausdorff distance. Then, $r_n\cdot \op{Cac}(\mathbf{G}_n)$
also converges to $\op{Kac}(\mathbf{E})$, in the
sense of the Gromov-Hausdorff distance.
\end{proposition}
 
 \begin{rek} \label{remgeodesic}
 {\rm The
 cactus $\op{Kac}(\mathbf{E})$ is well defined because $\mathbf{E}$ must be a geodesic space.
 The latter property can be derived from \cite[Theorem 7.5.1]{BBI01}, using the fact that 
 the graphs $r_n\cdot G_n$ can be approximated by geodesic spaces as explained 
 in the forthcoming proof.}
 \end{rek}
 
 \proof This is essentially a consequence of Proposition \ref{LipschitzGH}. 
 We start with some simple observations. Let $\mathbf{G}=(V,\mathcal{E},\rho)$ be
 a pointed graph. By considering 
 the union of a collection $(I_{\{u,v\}})_{\{u,v\}\in \mathcal{E}}$ of unit segments indexed by $\mathcal{E}$
 (such that this union is a metric graph in the sense of 
 \cite[Section 3.2.2]{BBI01}),
we can construct a pointed geodesic compact metric space $(\Lambda(\mathbf{G}),d_{\Lambda(\mathbf{G})},\tilde\rho)$, 
such that the graph $\mathbf{G}$ (viewed as a pointed metric space)
is embedded isometrically in $\Lambda(\mathbf{G})$, and the Gromov-Hausdorff distance
 between $\mathbf{G}$ and $\Lambda(\mathbf{G})$ is bounded above by $1$. 
 
 A moment's thought shows that $\op{Cac}(\mathbf{G})$ is also embedded isometrically in
 $\op{Kac}(\Lambda(\mathbf{G}))$, and the Gromov-Hausdorff distance
 between $\op{Cac}(\mathbf{G})$ and $\op{Kac}(\Lambda(\mathbf{G}))$ is still bounded above by $1$. 
 
 We apply these observations to the graphs $\mathbf{G}_n$. By scaling, we get that the
 Gromov-Hausdorff distance between the metric spaces $r_n\cdot \mathbf{G}_n$ and $r_n\cdot \Lambda(\mathbf{G}_n)$
 is bounded above by $r_n$, so that the sequence $r_n\cdot \Lambda(\mathbf{G}_n)$ also converges to $\mathbf{E}$
 in the sense of the Gromov-Hausdorff distance. From Proposition \ref{LipschitzGH}, we now get that
 $\op{Kac}(r_n\cdot \Lambda(\mathbf{G}_n))$ converges to $\op{Kac}(\mathbf{E})$. On the other hand,
 the Gromov-Hausdorff distance beween $\op{Kac}(r_n\cdot \Lambda(\mathbf{G}_n))=r_n\cdot \op{Kac}(\Lambda(\mathbf{G}_n))$
 and $r_n\cdot \op{Cac}(\mathbf{G}_n)$ is bounded above by $r_n$, so that the convergence of the
 proposition follows. 
 \endproof
 
 \begin{corollary}
 \label{cactus-tree}
 Let $\mathbf{E}$ be a pointed geodesic compact metric space. Then $\op{Kac}(\mathbf{E})$
 is a compact $\R$-tree.
 \end{corollary}
 
 \proof As a simple consequence of Proposition 7.5.5 in \cite{BBI01}, we can find a sequence $(r_n)_{n\geq 0}$
 of positive real numbers converging to $0$ and a sequence $(\mathbf{G}_n)_{n\geq 0}$ of pointed graphs,
 such that the rescaled graphs $r_n\cdot \mathbf{G}_n$ converge to $\mathbf{E}$ in the Gromov-Hausdorff sense.
 By Proposition \ref{convdiscactus}, $r_n\cdot \op{Cac}(\mathbf{G}_n)$ converges to $\op{Kac}(\mathbf{E})$
 in the Gromov-Hausdorff sense. Using the notation of the preceding proof, it also holds that
 $r_n\cdot\Lambda(\op{Cac}(\mathbf{G}_n))$ converges to $\op{Kac}(\mathbf{E})$. Proposition
 \ref{cakd} then implies that $r_n\cdot\Lambda(\op{Cac}(\mathbf{G}_n))$ is a (compact) $\R$-tree. The desired result
 follows since the set of all compact $\R$-trees is known to be closed for the Gromov-Hausdorff topology (see e.g. 
 \cite[Lemma 2.1]{EPW}).
 \endproof

 \subsection{Another approach to the continuous 
   cactus}\label{sec:furth-remarks-cont-1}

In this section, we present an alternative definition of the continuous cactus, which
 gives a different perspective on the previous results, and in particular on
 Corollary \ref{cactus-tree}. Let $\mathbf{E}=(E,d,\rho)$ be a pointed geodesic 
 compact metric space, and for $r\geq 0$, let
$$\mathbf{B}(r)=\{x\in E:d(\rho,x)<r\}\, ,\qquad \ov{\mathbf{B}}(r)=\{x\in E:d(\rho,x)\leq
r\}\, ,$$ be respectively  the open and the closed ball of radius $r$ centered at
$\rho$. We let $\op{Kac}'(\mathbf{E})$ be the set of all subsets of $E$
that are (non-empty) connected components of the closed set $\mathbf{B}(r)^c$,
for some $r\geq 0$ (here, $A^c$ denotes the complement of the set
$A$). Note that all elements of $\op{Kac}'(\mathbf{E})$ are themselves
closed subsets of $E$.

For every $C\in \op{Kac}'(\mathbf{E})$, we let
$$h(C)=d(\rho,C)=\inf\{d(\rho,x):x\in C\}\, .$$ 
Since $E$ is path-connected, $h(C)$ is also the unique real $r\geq 0$
such that $C$ is a connected component of $\mathbf{B}(r)^c$.

Note that $\op{Kac}'(\mathbf{E})$ is partially ordered by the
relation 
$$C\preceq C'\iff C'\subseteq C$$ and has a unique minimal element
$E$. Every totally ordered subset of $\op{Kac}'(\mathbf{E})$ has a
supremum, given by the intersection of all its elements. To see this,
observe that if $(C_i)_{i\in I}$ is a  totally ordered subset of $\op{Kac}'(\mathbf{E})$
then we can choose a sequence $(i_n)_{n\geq 1}$ taking values in $I$ such that
the sequence $(h(C_{i_n}))_{n\geq 1}$ is non-decreasing and converges to $r_{\op{max}}:=\sup\{h(C_i):i\in I\}$. Then
the intersection
$$\bigcap_{n= 1}^\infty C_{i_n}$$
is non-empty, closed and connected as the intersection of a decreasing sequence of 
non-empty closed connected sets in a compact space, and it easily follows that this intersection is
a connected component of $\mathbf{B}(r_{\op{max}})^c$ and coincides with the 
intersection of all $C_i$, $i\in I$.
At this point, it is crucial that elements of $\op{Kac}'(\mathbf{E})$ are
closed, and this 
is one of the reasons why one considers complements of {\em
  open} balls in the definition of $\op{Kac}'(\mathbf{E})$. 
  
  In
particular, for every $C,C'\in \op{Kac}'(\mathbf{E})$ , the infimum $C\wedge C'$
makes sense as the supremum of all $C''\in\op{Kac}'(\mathbf{E})$
such that $C''\preceq C$ and $C''\preceq C'$, and $h(C\wedge C')$ is the maximal 
value of $r$
such that $C$ and $C'$ are contained in the same connected component
of $\mathbf{B}(r)^c$. 

Moreover, if $C\in \op{Kac}'(\mathbf{E})$, the set $\{C'\in
\op{Kac}'(\mathbf{E}):C'\preceq C\}$ is isomorphic as an ordered set
to the segment $[0,h(C)]$, because for every $t\in [0,h(C)]$ there is
a unique $C'\in \op{Kac}'(\mathbf{E})$ with $h(C')=t$ and $C\subset
C'$. 

Finally, $h:\op{Kac}'(\mathbf{E})\to \R_+$ is an increasing function,
inducing a bijection from every segment of the partially ordered set
$\op{Kac}'(\mathbf{E})$ to a real segment. It follows from general
results (see Proposition 3.10 in \cite{FaJo04}) that the set
$\op{Kac}'(\mathbf{E})$ equipped with the distance
$$\op{d}^{\mathbf{E}}_{\op{Kac}'}(C,C')=h(C)+h(C')-2h(C\wedge C')$$
is an $\R$-tree rooted at $E=\mathbf{B}(0)^c$.
Note that $\op{d}^{\mathbf{E}}_{\op{Kac}'}(E,C)=h(C)$
for every $C\in \op{Kac}'(\mathbf{E})$.

\begin{proposition}\label{sec:furth-remarks-cont}
  The spaces $\op{Kac}'(\mathbf{E})$ and $\op{Kac}(\mathbf{E})$ are
  isometric pointed metric spaces.
\end{proposition}

\proof We consider the mapping from $E$ to $\op{Kac}'(\mathbf{E})$,
which maps $x$ to the connected component $C_x$ of
$\mathbf{B}(d(\rho,x))^c$ containing $x$. This mapping is clearly
onto: if $C\in \op{Kac}'(\mathbf{E})$, we have $C=C_x$ for any $x\in
C$ such that $d(\rho,x)=d(\rho,C)$. Let us show that this mapping is
an isometry from the pseudo-metric space
$(E,\op{d}^{\mathbf{E}}_{\op{Kac}})$ onto
$(\op{Kac}'(\mathbf{E}),\op{d}^{\mathbf{E}}_{\op{Kac}'})$.

Let $x,y\in E$ be given, and $\gamma:[0,1]\to E$ be a path from $x$ to
$y$. Let $t_0$ be such that $d(\rho,\gamma(t_0)) \leq
d(\rho,\gamma(t))$ for every $t\in [0,1]$. Then the path $\gamma$ lies
in a single path-connected component of
$\mathbf{B}(d(\rho,\gamma(t_0)))^c$, entailing that $x$ and $y$ are in the
same connected component of this set. Consequently, $h(C_x\wedge
C_y)\geq d(\rho,\gamma(t_0))$, and since obviously $h(C_x)=d(x,\rho)$,
$$\op{d}^{\mathbf{E}}_{\op{Kac}'}(C_x,C_y)\leq
d(\rho,x)+d(\rho,y)-2\inf_{t\in [0,1]}d(\rho,\gamma(t))\, .$$
Taking the infimum over all $\gamma$ gives 
\begin{equation}\label{eq:1}
\op{d}^{\mathbf{E}}_{\op{Kac}'}(C_x,C_y)\leq
\op{d}^{\mathbf{E}}_{\op{Kac}}(x,y)\, .
\end{equation}
Let us verify that the reverse inequality also holds. If $h(C_x\wedge
C_y)>0$ and $\eps\in (0,h(C_x\wedge C_y))$, the infimum
$C_x\wedge C_y$ is contained in some connected component of
$\ov{\mathbf{B}}(h(C_x\wedge C_y)-\eps)^c$. Since the latter set is open,
and $E$ is a geodesic space, hence locally path-connected, we
deduce that this connected component is in fact path-connected, and
since it contains $x$ and $y$, we can find a path $\gamma$ from $x$ to
$y$ that remains in $\ov{\mathbf{B}}(h(C_x\wedge C_y)-\eps)^c$. This
entails that 
$$\op{d}^{\mathbf{E}}_{\op{Kac}}(x,y)\leq
\op{d}^{\mathbf{E}}_{\op{Kac}'}(C_x,C_y)+\eps\, ,$$ and letting
$\eps\to 0$ yields the bound $\op{d}^{\mathbf{E}}_{\op{Kac}'}(C_x,C_y)\geq
\op{d}^{\mathbf{E}}_{\op{Kac}}(x,y)$. The latter bound remains true
when $h(C_x\wedge C_y)=0$, since in that case 
$C_x\wedge C_y=E$ and 
$\op{d}^{\mathbf{E}}_{\op{Kac}'}(C_x,C_y)=h(C_x)+h(C_y)=d(\rho,x)+d(\rho,y)$.

From the preceding observations, we directly obtain that $x\mapsto C_x$ induces a quotient
mapping from $\op{Kac}(\mathbf{E})$ onto $\op{Kac}'(\mathbf{E})$,
which is an isometry and maps (the class of) $\rho$ to
$E$. \endproof

 \begin{figure}[!h]
 \label{counterex}
 \begin{center}
\includegraphics[width=10cm,height=8cm]{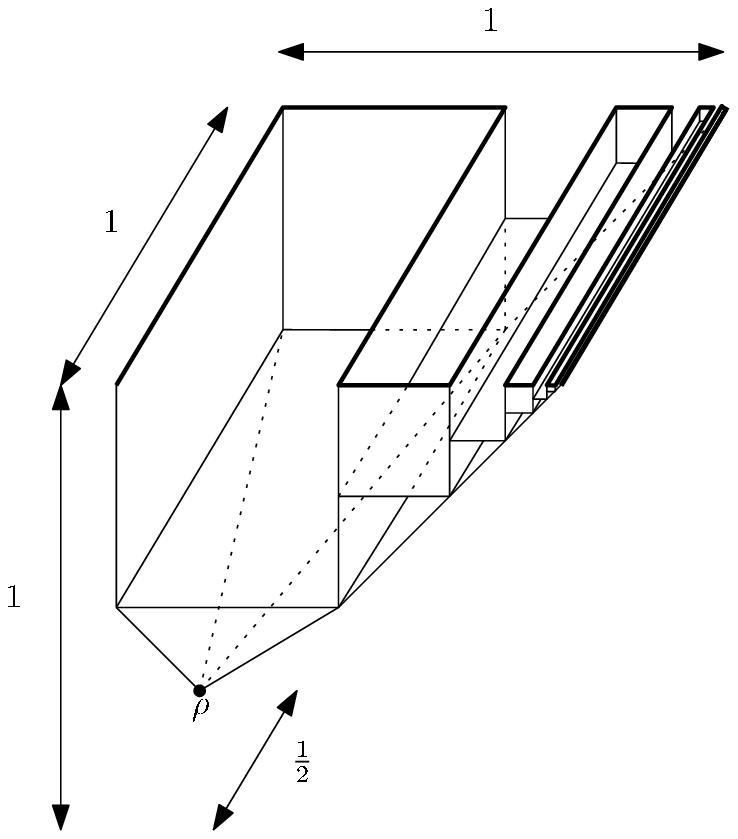}
\end{center}
\caption{An example of a geodesic compact metric space $E$, such that the complement
of the open ball of
radius $1$ centered at the distinguished point $\rho$ is connected but
not path-connected. Here $E$ is a compact subset of $\R^3$ and is equipped with
the intrinsic distance associated with the $L^\infty$-metric $\delta((x_1,x_2,x_3),(y_1,y_2,y_3))
=\sup\{|x_i-y_i|,i=1,2,3\}$. For this distance, the sphere of radius $1$ centered at $\rho$, which coincides
with the complement
of the open ball of
radius $1$,
consists of the union of the bold lines at the top of the figure.
}
\end{figure}

\begin{rek} {\rm The discrete cactus of a graph can be defined in an analogous way as
above, using the notion of graph connectedness instead of
connectedness in metric spaces. }
\end{rek}

Let us return to Remark \ref{minimiz} about the existence, for given $x,y\in E$, of a minimizing path 
$\gamma:[0,1]\to E$ going from $x$ to $y$, such that
$$\op{d}^{\mathbf{E}}_{\op{Kac}}(x,y)=d(\rho,x)+d(\rho,y)-2\min_{0\leq
  t\leq 1}d(\rho,\gamma(t)).$$ 
  With the notation of the
previous proof, it may happen that the closed set $C_x\wedge C_y$ is
connected without being path-connected: Fig.2 suggests an example 
of this phenomenon. In that event, if $x$ and $y$ cannot be connected
by a continuous path that stays in $C_x\wedge C_y$, there exists no
minimizing path.

 \section{The Brownian cactus}
 \label{Brcactus}
 
 In this section, we define the Brownian cactus and we show that it is the continuous
 cactus associated with the (random) compact metric space called
 the Brownian map. The Brownian map has been studied in \cite{IM} as the 
 limit in distribution, along suitable sequences, of rescaled
 $2p$-angulations chosen uniformly at random. 
 We first recall  some basic facts about the Brownian map.
 
 We let $\be=(\be_t)_{0\leq t\leq 1}$ be a Brownian excursion with
 duration $1$. For our purposes it is crucial to view $\be$ as the
 coding function for the random continuous tree known as the
 CRT. Precisely, we define a pseudo-distance
 $\op{d}_\be$ on $[0,1]$ by setting for every $s,t\in[0,1]$,
 $$\op{d}_\be(s,t) =\be_s+\be_t - 2\,\min_{s\wedge t\leq r\leq s\vee t} \be_r$$
 and we put $s\sim_\be t$ iff $\op{d}_\be(s,t)=0$. The CRT is 
 defined as the quotient metric space $\T_\be:=[0,1]\,/\!\sim_\be$,
and is equipped with the induced metric $\op{d}_\be$. 
Then $(\T_\be,\op{d}_\be)$ is a random (compact) $\R$-tree. We write
 $p_\be:[0,1]\la \T_\be$ for the canonical projection, and we define the mass
 measure (or volume measure) $\op{Vol}$ on the CRT as the image 
 of Lebesgue measure on $[0,1]$ under $p_\be$. 
 For every $a,b\in\T_\be$,
 we let $\llbracket a,b \rrbracket$ be the range of the
 geodesic path from $a$ to $b$ in $\T_\be$: This is the line
 segment between $a$ and $b$ in the tree $\T_\be$. We will need the following simple
 fact, which is easily checked from the definition of 
 $\op{d}_\be$. Let $a,b\in \T_\be$, and let $s,t\in[0,1]$ be such that
 $p_\be(s)=a$ and $p_\be(t)=b$. Assume for definiteness that 
 $s\leq t$. Then $\llbracket a,b\rrbracket$ exactly consists of the
 points $c$ that can be written as $c=p_\be(r)$, with $r\in[s,t]$
 satisfying
 $$\be_r=\max\Big(\min_{u\in[s,r]} \be_u,\min_{u\in[r,t]} \be _u\Big).$$

 Conditionally given $\be$, we introduce the
 centered Gaussian process $(Z_t)_{0\leq t\leq 1}$ with continuous sample paths
 such that
 $${\rm cov}(Z_s,Z_t)= \min_{s\wedge t\leq r\leq s\vee t} \be_r.$$
 It is easy to verify that a.s. for every $s,t\in[0,1]$ the condition 
 $s\sim_\be t$ implies that $Z_s=Z_t$. Therefore we may and will
 view $Z$ as indexed by the CRT $\T_\be$. In fact, it is natural to
 interpret $Z$ as Brownian motion indexed by the CRT.
 We will write indifferently $Z_a=Z_t$ when $a\in\T_\be$ and 
 $t\in[0,1]$ are such that $a=p_\be(t)$. 
 
 We set
 $$\un Z:= \min_{t\in[0,1]} Z_t.$$
 One can then prove \cite{MM,LGW} that a.s. there exists a unique
 $s_*\in[0,1]$ such that $Z_{s_*}=\un Z$. We put
 $$a_*=p_\be(s_*).$$
 
We now define an equivalence relation on the CRT. For every
$a,b\in\T_\be$, we put $a\approx b$ if and only if 
there exist $s,t\in[0,1]$ such that $p_\be(s)=a$, $p_\be(t)=b$, and
$$Z_r\geq Z_s=Z_t\;,\qquad\hbox{for every }r\in[s,t].$$
Here and later we make the convention that when 
$s>t$, the notation $r\in[s,t]$ means $r\in [s,1]\cup [0,t]$.

It is not obvious that $\approx$ is an 
equivalence relation. This follows from Lemma 3.2 in \cite{LGP}, which
shows that with probability one, for every
distinct $a,b\in\T_\be$, the property $a\approx b$ may only
hold if $a$ and $b$ are leaves of $\T_\be$, and then $p_\be^{-1}(a)$
and $p_\be^{-1}(b)$ are both singletons.

The Brownian map is now defined as the quotient space
$$m_\infty := \T_\be\, / \!\approx$$
which is equipped with the quotient topology. We 
write $\Pi:\T_\be \la m_\infty$ for the canonical projection, and
we put $\rho_*=\Pi(a_*)$. We also let $\lambda$
be the image of $\op{Vol}$ under $\Pi$, and we interpret
$\lambda$ as the volume measure on $m_\infty$. For every $x\in m_\infty$,
we set $Z_x=Z_a$, where $a\in\T_\be$ is such that $\Pi(a)=x$ (this definition does
not depend on the choice of $a$).

A key result of \cite{IM} states the Brownian map, equipped 
with a suitable metric $D$, appears as the limit in
distribution of rescaled random $2p$-angulations. More precisely, let
$p\geq 2$ be an integer, and for every $n\geq 1$, let $m_n$
be uniformly distributed over the class of all 
rooted $2p$-angulations with $n$ faces. Write $V(m_n)$ for
the vertex set of $m_n$, which is equipped with the
graph distance $\op{d}_{\op{gr}}^{m_n}$, and let $\rho_n$
denote the root vertex of $m_n$. Then, from any strictly increasing sequence
of positive integers we can
extract a suitable subsequence $(n_k)_{k\geq 1}$ such that the following convergence holds in distribution in
the Gromov-Hausdorff sense,
\begin{equation}
\label{conv-to-Bmap}
\Big(V(m_{n_k}),\Big(\frac{9}{4p(p-1)}\Big)^{1/4}(n_k)^{-1/4}\op{d}_{\op{gr}}^{m_{n_k}},\rho_{n_k}\Big)
 \build{\la}_{k\to\infty}^{\rm (d)} (m_\infty, D,\rho_*)
 \end{equation}
 where $D$ is a metric on the space $m_\infty$
 that satisfies the following properties:
 \begin{enumerate}
 \item For every $a\in\T_\be$,
 $$D(\rho_*,\Pi(a))= Z_a - \un Z.$$
 \item For every $a,b\in \T_\be$ and every $s,t\in[0,1]$ such that $p_\be(s)=a$ and $p_\be(t)=b$,
 $$D(\Pi(a),\Pi(b))\leq Z_s + Z_t -2 \min_{r\in[s,t]} Z_r\;.$$
 \item For every $a,b\in \T_\be$,
 $$D(\Pi(a),\Pi(b)) \geq Z_a + Z_b - 2 \min_{c\in\llbracket a,b \rrbracket} Z_c\;.$$
 \end{enumerate}
 
 Notice that in Property 2 we make the same convention as above for the notation
 $r\in[s,t]$ when $s>t$. The preceding statements can be found in Section 3 of \cite{IM}
 (see in particular \cite[Theorem 3.4]{IM}), with the exception of Property 3. We refer to Corollary 3.2 in \cite{Geo} for
the latter property. By the argument in Remark \ref{remgeodesic}, the metric  space $(m_\infty, D)$ is a geodesic space a.s.
 
 The limiting metric $D$ in (\ref{conv-to-Bmap}) may depend on the integer $p$ and on
 the choice of the subsequence $(n_k)$. However, we will see that
the cactus of the Brownian map is well defined independently of $p$ and of the
chosen  subsequence, and in fact coincides with the Brownian
 cactus that we now introduce.
 
 \begin{definition}
 \label{Bcactus}
 The Brownian cactus  $\op{KAC}$ is the random metric space 
 defined as the quotient space of $\T_\be$
 for the equivalence relation 
 $$a\asymp b\quad\hbox{iff}\quad Z_a=Z_b=\min_{c\in\llbracket a,b \rrbracket} Z_c$$
 and equipped with the distance induced by
 $$\op{d}_{\op{KAC}}(a,b)= Z_a + Z_b - 2 \min_{c\in\llbracket a,b \rrbracket} Z_c\;,\ \hbox{for every }a,b\in\T_\be.$$
 We view $\op{KAC}$ as a pointed metric space whose root is the equivalence class of $a_*$.
 \end{definition}
 
 It is an easy matter to verify that $\op{d}_{\op{KAC}}$ is a pseudo-distance on $\T_\be$, and that 
 $\asymp$ is the associated equivalence relation. 
 
 We write $\bm_\infty$ for the pointed metric space $(m_\infty,D,\rho_*)$ appearing in (\ref{conv-to-Bmap}).
 
 \begin{proposition}
 \label{identKactus}
Almost surely, $\op{Kac}(\bm_\infty)$ is isometric to $\op{KAC}$.
 \end{proposition}
 
 \proof We first need to identify the pseudo-distance $\op{d}_{\op{Kac}}^{\bm_\infty}$
 (see subsection \ref{sec:contiK}). Let $x,y\in m_\infty$ and choose 
 $a,b\in\T_\be$ such that $x=p_\be(a)$ and $y=p_\be(b)$. 
 If $\gamma:[0,1]\la m_\infty$ is a continuous path 
 such that $\gamma(0)=x$ and $\gamma(1)=y$, Proposition 3.1
 in \cite{Geo} ensures that
 $$\min_{0\leq t\leq 1} Z_{\gamma(t)} \leq \min_{c\in\llbracket a,b\rrbracket} Z_c.$$
 Using Property 1 above, it follows that
 $$\min_{0\leq t\leq 1} D(\rho_*,\gamma(t)) \leq \min_{c\in\llbracket a,b\rrbracket} (Z_c-\un Z).$$
 Since this holds for any continuous curve $\gamma$ from $x$ to $y$ 
 in $m_\infty$, we get from the definition of $\op{d}_{\op{Kac}}^{\bm_\infty}$ that
 $$\op{d}_{\op{Kac}}^{\bm_\infty}(x,y)\geq (Z_a-\un Z)+(Z_b-\un Z) - 2 \min_{c\in\llbracket a,b\rrbracket} (Z_c-\un Z)= \op{d}_{\op{KAC}}(a,b).$$
 The corresponding upper bound is immediately obtained by letting $\gamma$
 be the image under $\Pi$ of the (rescaled) geodesic path from $a$ to $b$ in the tree $\T_{\be}$. Note that
 the resulting path from $x$ to $y$ in $m_\infty$ is continuous because the projection $\Pi$
 is so. Summarizing, we have obtained that, for every $a,b\in\T_{\be}$,
 \begin{equation}
 \label{identitKac}
 \op{d}_{\op{Kac}}^{\bm_\infty}(\Pi(a),\Pi(b))= \op{d}_{\op{KAC}}(a,b).
 \end{equation}
 In particular, the property $a\asymp b$ holds if and only if $\Pi(a)\build{\asymp}_{}^{\bm_\infty} \Pi(b)$. 
Hence, the composition of the canonical projections from $\T_\be$ onto $m_\infty$
 and from $m_\infty$ onto $\op{Kac}(\bm_\infty)$ induces a one to-one mapping from
 $\op{KAC}=\T_\be/\asymp$ onto $\op{Kac}(\bm_\infty)$. By (\ref{identitKac}) this mapping
 is an isometry, which completes the proof. \endproof
 
 \medskip
 Recall the notation $m_n$ for a random planar map uniformly distributed over the 
 set of all rooted $2p$-angulations with $n$ faces, and $\rho_n$ for the root vertex of $m_n$. As explained at the end of subsection
 \ref{sec:disccac}, we can associate a pointed graph with $m_n$, such that the distinguished point of this graph 
 is $\rho_n$. We write $\bm_n$ for this pointed graph. 
 \begin{corollary}
 \label{conv-to-cactus}
 We have 
 $$\Big(\frac{9}{4p(p-1)}\Big)^{1/4}n^{-1/4}\cdot \op{Cac}(\bm_n)
 \build{\la}_{n\to\infty}^{\rm (d)} \op{KAC}$$
 in the Gromov-Hausdorff sense.
 \end{corollary}
 
 In contrast with (\ref{conv-to-Bmap}), the convergence of the corollary does not
 require the extraction of a subsequence.
 
 \medskip
 \proof It is sufficient to prove that, from any strictly increasing sequence 
 of positive integers we can extract a subsequence $(n_k)$ such that the desired 
 convergence holds along this subsequence. To this end,
 we extract the subsequence $(n_k)$ so that (\ref{conv-to-Bmap}) holds. 
 By Proposition \ref{convdiscactus}, we have 
 then
$$ \Big(\frac{9}{4p(p-1)}\Big)^{1/4}(n_k)^{-1/4}\cdot \op{Cac}(\bm_{n_k})
\build{\la}_{k\to\infty}^{\rm (d)} \op{Kac}(\bm_\infty).$$
By Proposition \ref{identKactus}, the limiting distribution is that of $\op{KAC}$,
independently of the subsequence that we have chosen. This completes
the proof. \endproof
 
 \medskip
 In the next section, we will see that the convergence of the corollary holds
 for much more general random planar maps.
 
 \section{Convergence of cactuses associated with random planar maps}
 \label{convcactusmaps}
 
 \subsection{Planar maps and bijections with trees}
 \label{biject}
 
 We denote 
 the set of all rooted and pointed planar maps by $\m_{r,p}$. 
 As in \cite{MieInvar}, it is convenient for technical reasons to make
 the convention that
 $\m_{r,p}$ contains the ``vertex map'', denoted by $\dagger$, which
 has no edge and only one vertex ``bounding'' a face of degree $0$.
With the exception of $\dagger$, a  
planar map in $\m_{r,p}$ has
 at least one edge. An element of $\m_{r,p}$ other than $\dagger$
 consists of a planar map $m$ together with an oriented edge $e$ (the root edge) and 
 a distinguished vertex $\rho$. We write $e_-$ and $e_+$ for the origin and 
 the target of the root edge $e$. Note that we may have $e_-=e_+$ if
 $e$ is a loop. 
 
 As previously, we denote the graph distance on the vertex set $V(m)$ of $m$ by $\op{d}_{\op{gr}}^m$.
 We say that the rooted and pointed planar map $(m,e,\rho)$ is positive, respectively negative, respectively null
 if $\op{d}_{\op{gr}}^m(\rho,e_+)=\op{d}_{\op{gr}}^m(\rho,e_-) +1$, resp. $\op{d}_{\op{gr}}^m(\rho,e_+)=\op{d}_{\op{gr}}^m(\rho,e_-) -1$, resp.
 $\op{d}_{\op{gr}}^m(\rho,e_+)=\op{d}_{\op{gr}}^m(\rho,e_-)$. We make the convention that the
 vertex map $\dagger$ is positive. We write $\m_{r,p}^+$, resp. $\m_{r,p}^-$, resp. $\m_{r,p}^0$
 for the set of all positive, resp. negative, resp. null, rooted and pointed planar maps.
 Reversing the orientation of the root edge yields an obvious bijection between the sets
 $\m_{r,p}^+$ and $\m_{r,p}^-$, and for this reason we will mainly discuss $\m_{r,p}^+$
 and $\m_{r,p}^0$ in what follows.
 
 We will make use of the Bouttier-Di Francesco-Guitter bijection \cite{BDG}
 between $\m_{r,p}^+\cup \m_{r,p}^0$ and a certain set of multitype labeled trees called 
 mobiles. In order to describe this bijection, we use the standard formalism for plane trees, as found in 
 Section 1.1 of \cite{tree} for instance. In this formalism, vertices are elements of
 the set
 $$\mathcal{U}=\bigcup_{n=0}^\infty \N^n$$
 of all finite sequences of positive integers, including the empty sequence $\varnothing$
 that serves as the root vertex of the tree. 
A plane tree $\tau$ is a finite subset of $\mathcal{U}$ that satisfies the following three conditions:
\begin{enumerate}
\item $\varnothing\in\tau$.
\item For every $u=(i_1,\ldots,i_k)\in\tau\setminus\{\varnothing\}$, the sequence $(i_1,\ldots,i_{k-1})$ (the ``parent'' of $u$)
also belongs to $\tau$.
\item For every $u=(i_1,\ldots,i_k)\in\tau$, there
exists an integer $k_u(\tau)\geq 0$ (the ``number of children'' of $u$) such that
the vertex $(i_1,\ldots,i_k,j)$ belongs to $\tau$ if
and only if $1\leq j\leq k_u(\tau)$. 
\end{enumerate}
The generation of $u=(i_1,\ldots,i_k)$ is denoted by $|u|=k$. 
The notions of an ancestor and a descendant in the tree $\tau$
are defined in an obvious way.

 We will be interested in four-type plane trees,
 meaning that each vertex is assigned a type which can be $1,2,3$ or $4$. 
 
  We next introduce mobiles following the presentation in \cite{MieInvar}, with 
  a few minor modifications.
 We consider a four-type plane tree $\tau$ satisfying the following properties:
\begin{enumerate}
\item[(i)] The root vertex $\varnothing$ is of type $1$ or of type $2$.
\item[(ii)] The children of any vertex of type $1$ are of type $3$.
\item[(iii)] Each individual of type $2$ and which is not the root vertex of the tree has exactly one child of type $4$
and no other child.
If the root vertex is of type $2$, it has exactly two children, both of type $4$.
\item[(iv)] The children of individuals of type $3$ or $4$ can only be of type $1$ or $2$.
\end{enumerate}
 Let $\tau_{(1,2)}$ be the set of all vertices of $\tau$ at even generation
(these are exactly the vertices of type $1$ or $2$). An admissible labeling
of $\tau$ is a collection of integer labels $(\ell_u)_{u\in\tau_{(1,2)}}$
assigned to the vertices of type $1$ or $2$, such that the 
following properties hold:
 \begin{enumerate}
\item[a.] $\ell_\varnothing=0$.
\item[b.] Let $u$ be a vertex of type $3$ or $4$, let $u_{(1)},\ldots,u_{(k)}$
be the children of $u$ (in lexicographical order) and let $u_{(0)}$ be
the parent of $u$. Then, for every $i=0,1,\ldots,k$,
$$\ell_{u_{(i+1)}}\geq \ell_{u_{(i)}}-1$$
with the convention $u_{(k+1)}=u_{(0)}$. Moreover, for every $i=0,1,\ldots,k$ 
such that $u_{(i+1)}$ is of type $2$, we have
$$\ell_{u_{(i+1)}}\geq \ell_{u_{(i)}}.$$
\end{enumerate}

By definition, a \emph{mobile} is a pair $(\tau, (\ell_{u})_{u\in \tau_{(1,2)}})$ consisting 
of a four-type plane tree satisfying the preceding conditions (i)--(iv), and 
an admissible labeling of $\tau$. We let $\TT_+$ be the set of all
mobiles such that the root vertex of $\tau$ is of type $1$. We also let
$\TT_0$ be the set of all mobiles such that the root vertex is of type $2$.

\begin{rek} {\rm Our definition of admissible labelings is slightly different from the ones
that are used in \cite{MieInvar} or \cite{MW}. To recover the definitions of 
\cite{MieInvar} or \cite{MW}, just subtract $1$ from the label of each vertex of type $2$. 
Because of this difference, our construction of the bijections between maps and 
trees will look slightly different from the ones in \cite{MieInvar} or \cite{MW}.}
\end{rek}

The Bouttier-Di Francesco-Guitter construction provides bijections
between the set $\TT_+$ and the set $\m^+_{r,p}$ on one hand, 
between the set $\TT_0$ and the set $\m^0_{r,p}$ on the other hand.
Let us describe this construction in the first case. 

We start from a 
mobile $(\tau, (\ell_{u})_{u\in \tau_{(1,2)}})\in \TT_+$. In the case when 
$\tau=\{\varnothing\}$, we decide by convention that the associated 
planar map is the vertex map $\dagger$. Otherwise,
let $p\geq 1$
be the number of edges of $\tau$ ($p=\#\tau -1$). The contour sequence of 
$\tau$ is the sequence $v_0,v_1,\ldots,v_{2p}$ of vertices
of $\tau$ defined inductively as follows. First $v_0=\varnothing$.
Then, for every $i\in\{0,1,\ldots,2p-1\}$, $v_{i+1}$ is either the first child
of $v_i$ that has not yet appeared among $v_0,v_1,\ldots,v_i$, or
if there is no such child, the parent of $v_i$. It is easy to see that
this definition makes sense and $v_{2p}=\varnothing$. Moreover all
vertices of $\tau$ appear in the sequence $v_0,v_1,\ldots,v_{2p}$,
and more precisely the number of occurences of a vertex $u$
of $\tau$ is equal to the multiplicity of $u$ in $\tau$. In fact, each index $i$
such that $v_i=u$ corresponds to one corner of the vertex
$u$ in the tree $\tau$ : We will abusively call it the corner $v_i$. We also
introduce the modified contour sequence of $\tau$ 
as the sequence $u_0,u_1,\ldots,u_p$ defined by
$$u_i=v_{2i}\;,\qquad \forall i=0,1,\ldots,p.$$
By construction, the vertices appearing in the modified contour
sequence are exactly the vertices of $\tau_{(1,2)}$. We extend the
modified contour sequence periodically by setting 
$u_{p+i}=u_i$ for $i=1,\ldots,p$. Note that the 
properties of labels entail $\ell_{u_{i+1}}\geq \ell_{u_i} -1$ for $i=0,1,\ldots,2p-1$.

\begin{figure}[!h]
\begin{center}
\includegraphics[width=11cm,height=10cm]{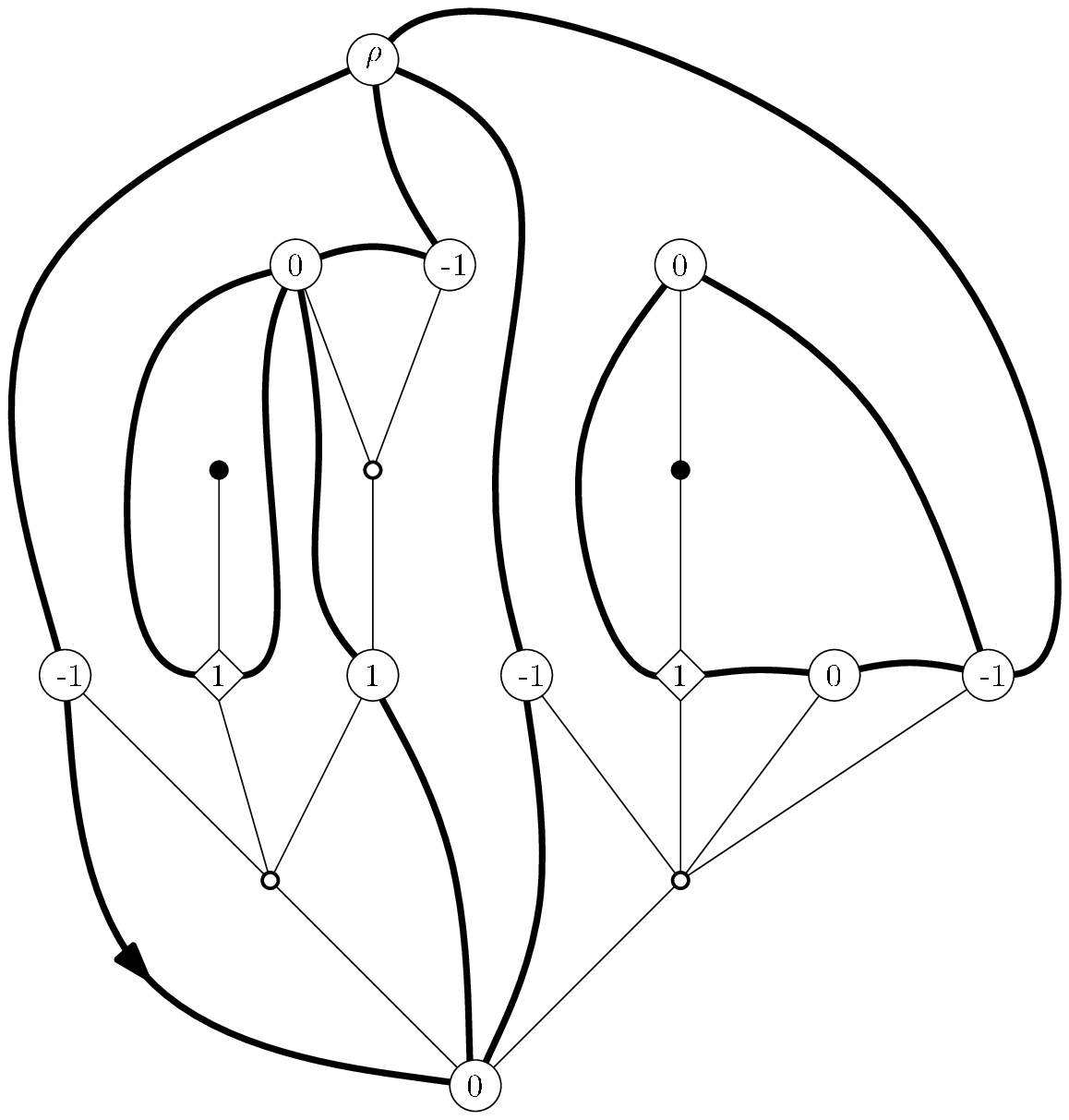}
\end{center}
\caption{A mobile $(\tau,(\ell_u)_{u\in\tau_{(1,2)}})$ in $\TT_+$ and its image $m$ under the BDG bijection.
Vertices of type $1$ are represented by big circles, vertices of type $2$ by lozenges, vertices of type $3$ by small circles
and vertices of type $4$ by small black disks. The edges of the tree $\tau$ are represented by thin lines, and the edges 
of the planar map $m$  by thick curves. In order to get the planar map $m$ one needs to erase the vertices of type $2$
and, for each of these vertices, to merge its two incident edges into a single one. The root edge is at the bottom left.
}
\end{figure}

To construct the edges of the rooted and pointed planar map $(m,e,\rho)$
associated with the mobile $(\tau, (\ell_{u})_{u\in \tau_{(1,2)}})\in \TT_+$
we proceed as follows. We first embed the tree $\tau$ in the plane
in a way consistent with the planar order. We then add an extra vertex 
of type $1$, which we call $\rho$. Then, for every $i=0,1,\ldots,p-1$:
\begin{enumerate}
\item[(i)] If 
$$\ell_{u_i}=\min_{0\leq k\leq p} \ell_{u_k}$$
 we draw an edge between the corner $u_i$ and $\rho$.
 \item[(ii)] If
 $$\ell_{u_i}>\min_{0\leq k\leq p} \ell_{u_k}$$
 we draw an
edge between the corner $u_i$ and the corner $u_j$, where
$j=\min\{k\in\{i+1,\ldots,i+p-1\}: \ell_{u_k}=\ell_{u_i}-1\}$. Because of
property b. of the labeling, the vertex $u_j$ must be of type $1$.
\end{enumerate}
The construction can be made in such a way that edges do not intersect,
and do not intersect the edges of the tree $\tau$. Furthermore each face of the resulting planar
map contains exactly one vertex of type $3$ or $4$, and both the parent and the
children of this vertex are incident to this face. See Fig.2 for an example.

The resulting planar map is bipartite with vertices either of type $1$
or of type $2$. Furthermore, the fact that in the tree $\tau$ each vertex
of type $2$ has exactly one child, and the labeling rules imply 
that each vertex of type $2$ is incident to exactly two edges of the map,
which connect it to two vertices of type $1$, which may be the same (these vertices
of type $1$ will be said to be associated with the vertex of type $2$ we are considering).
Each of these edges corresponds in the preceding construction to one
of the two corners of the vertex of type $2$ that we consider. To complete the
construction, we just erase all vertices of type $2$ and for each of these
we merge its two incident edges into a single edge connecting the two associated vertices
of type $1$. In this way we get a (non-bipartite in general)
planar map $m$. Finally we decide that the root edge $e$ of the map is the first edge drawn
in the construction, oriented in such a way that $e_+=\varnothing$, and we let  the distinguished vertex
of the map be the vertex $\rho$.
Note that vertices of the map $m$ that are different from the distinguished vertex $\rho$
are exactly  the vertices of type $1$ in the tree $\tau$. In other words,
the vertex set $V(m)$ is identified with the set $\tau_{(1)}\cup\{\rho\}$, where
$\tau_{(1)}$ denotes the set of all vertices of $\tau$ of type $1$. 

The mapping $(\tau, (\ell_{u})_{u\in \tau_{(1,2)}})\la (m,e,\rho)$ that we have just described
is indeed a bijection from $\TT_+$ onto $\m_{r,p}^+$. We can construct a similar bijection
from $\TT_+$ onto $\m_{r,p}^-$ by the same construction, with the minor modification
that we orient the root edge in such a way that $e_-=\varnothing$.

Furthermore we can also adapt the preceding construction in order to get 
a bijection from $\TT_0$ onto $\m_{r,p}^0$. The construction of edges 
of the map proceeds in the same way, but the root edge is now obtained
as the edge resulting of the merging of the two edges incident to
$\varnothing$ (recall that for a tree in $\TT_0$ the root $\varnothing$
is a vertex of type $2$ that has exactly two children, hence also two corners). 
The orientation of the root edge is chosen according to some
convention: For instance, one may decide that the ``half-edge'' coming from
the first corner of $\varnothing$ corresponds to the origin of the
root edge.

In all three cases, distances in the planar map $m$ satisfy the following
key property: For every vertex $u\in\tau_{(1)}$, we have
\begin{equation}
\label{distance-label}
\op{d}_{\op{gr}}^m(\rho,u)= \ell_u - \min \ell +1
\end{equation}
where $\min \ell$ denotes the minimal label on the tree $\tau$. In the left-hand side
$u$ is viewed as a vertex of the map $m$, in agreement with the
preceding construction.

The three bijections we have described are called the BDG bijections. In the
remaining part of this section, we fix a mobile $(\tau,(\ell_{u})_{u\in\tau_{(1,2)}})$ 
belonging to $\TT_+$ (or to $\TT_0$) and its
image $(m,e,\rho)$ under the relevant BDG bijection.

\begin{rek}
{\rm We could have defined the BDG bijections without distinguishing between types $3$ and $4$. 
However, this distinction will be important in the next section when we consider 
random planar maps and the associated (random) trees. We will see that these random
trees are Galton-Watson trees with a different
offspring distribution for vertices of type $3$ than for vertices of
type $4$. }
\end{rek}

If $u,v\in \tau_{(1,2)}$, we denote by $\llbracket u,v\rrbracket$
the set of all vertices of type $1$ or $2$ that lie on the geodesic path from $u$
to $v$ in the tree $\tau$. 

\begin{proposition}
\label{cactusLB}
For every $u,v\in V(m)\backslash \{\rho\}=\tau_{(1)}$, and every path 
$\gamma=(\gamma(0),\gamma(1),\ldots,\gamma(k))$ in $m$
such that $\gamma(0)=u$ and $\gamma(k)=v$, we have
$$\min_{0\leq i\leq k} \op{d}_{\op{gr}}^m(\rho,\gamma(i))
\leq \min_{w\in\llbracket u,v\rrbracket} \ell_w - \min \ell +1.$$
\end{proposition}

\proof We may assume that the path $\gamma$ does not visit $\rho$, since
otherwise the result is trivial. Using (\ref{distance-label}), the statement reduces to
$$\min_{0\leq i\leq k} \ell_{\gamma(i)}
\leq \min_{w\in\llbracket u,v\rrbracket} \ell_w.$$
So we fix $w\in\llbracket u,v\rrbracket$ and we verify that
$\ell_{\gamma(i)}\leq \ell_w$ for some $i\in\{0,1,\ldots,k\}$.
We may assume that $w\not = u$ and $w\not =v$. The removal 
of the vertex $w$ (and of the edges incident to $w$) disconnects the tree $\tau$ in several connected components.
Write $C$ for the connected component containing $v$, and note that
this component does not contain $u$. Then let $j\geq 1$ be the first integer
such that $\gamma(j)$ belongs to $C$. Thus $\gamma(j-1)\notin C$, $\gamma(j)\in C$
and the vertices $\gamma(j-1)$ and $\gamma(j)$ are
linked by an edge of the map $m$. From (\ref{distance-label}), we have $|\ell_{\gamma(j)}-\ell_{\gamma(j-1)}|\leq 1$.  
Now we use the fact that the edge between $\gamma(j-1)$ and $\gamma(j)$ is produced by
the BDG bijection. Suppose first that $\gamma(j-1)$ and $\gamma(j)$ have a
different label. In that case, noting that the modified contour sequence must visit
$w$ between any visit of $\gamma(j-1)$ and any visit of $\gamma(j)$, we easily
get that $\min\{\ell_{\gamma(j)},\ell_{\gamma(j-1)}\}\leq \ell_w$ (otherwise our
construction could not produce an edge from $\gamma(j-1)$ to $\gamma(j)$). 
A similar argument applies to the case when $\gamma(j-1)$ and $\gamma(j)$ have the
same label. In that case, the edge between $\gamma(j-1)$ and $\gamma(j)$ must
come from the merging of two edges originating from a vertex of $\tau$
of type $2$. This vertex of type $2$ has to belong to the set
$\llbracket \gamma(j-1),\gamma(j)\rrbracket$ (which contains $w$), because otherwise the
two associated vertices of type $1$ could not be $\gamma(j-1)$ and $\gamma(j)$.
It again follows from our construction that we must have 
$\min\{\ell_{\gamma(j)},\ell_{\gamma(j-1)}\}\leq \ell_w$. This completes the proof. \endproof

\medskip
In the next corollary, we write $\bm$ for the graph associated with the map $m$ (in the sense of subsection \ref{sec:disccac}), which is pointed at the distinguished 
vertex $\rho$. The notation $\op{d}^\bm_{\op{Cac}}$
then refers to the cactus distance for this pointed graph.

\begin{corollary}
\label{bounddiscretecactus}
Suppose that the degree of all faces of $m$ is bounded above by $D\geq 1$. Then,
for every $u,v\in V(m)\backslash \{\rho\}$, we have
$$\Big|\op{d}^\bm_{\op{Cac}}(u,v) -
\Big(\ell_u + \ell_v -2 \min_{w\in\llbracket u,v\rrbracket} \ell_w\Big)\Big| \leq 2D + 2.$$
\end{corollary}

\proof From the definition of the cactus distance $\op{d}^\bm_{\op{Cac}}$ and the preceding proposition,
we immediately get the lower bound
\begin{align*}
\op{d}^\bm_{\op{Cac}}(u,v)& \geq \op{d}_{\op{gr}}^m(\rho,u) + \op{d}_{\op{gr}}^m(\rho,v) - 2\Big(\min_{w\in\llbracket u,v\rrbracket} \ell_w - \min \ell +1\Big)\\
&= \ell_u + \ell_v -2 \min_{w\in\llbracket u,v\rrbracket} \ell_w,
\end{align*}
by (\ref{distance-label}).
In order to get a corresponding upper bound, let $\eta(0)=u,\eta(1),\ldots,\eta(k)=v$ be the vertices of type $1$
or $2$ belonging to the geodesic
path from $u$ to $v$ in the tree $\tau$, enumerated in their order of appearance on this path. Put $\tilde \eta(i)=\eta(i)$
if $\eta(i)$ is of type $1$, and if $\eta(i)$ is of type $2$, let $\tilde\eta(i)$ be one of the
two (possibly equal) vertices of type $1$ that are associated with $\eta(i)$ in the 
BDG bijection. Then the properties of the BDG bijection ensure that, for every
$i=0,1,\ldots,k-1$, the two vertices $\eta(i)$ and $\eta(i+1)$ lie on the boundary of
the same face of $m$ (the point is that, in the BDG construction, edges of the map $m$
are drawn in such a way that they do not cross edges of the tree $\tau$). 
From our assumption we have thus $\op{d}_{\op{gr}}^m(\tilde\eta(i),\tilde\eta(i+1))\leq D$
for every $i=0,1,\ldots,k-1$. Hence, we can find a path $\gamma$ in $m$ starting from $u$
and ending at $v$, such that
\begin{align*}
\min_j \op{d}_{\op{gr}}^m(\rho,\gamma(j))&\geq \min_{0\leq i\leq k} \op{d}_{\op{gr}}^m(\rho,\tilde\eta(i)) -D\\
&= \min_{0\leq i\leq k} \ell_{\tilde\eta(i)} - \min\ell +1 -D\\
&\geq \min_{0\leq i\leq k} \ell_{\eta(i)} - \min\ell -D\;.
\end{align*}
It follows that
\begin{align*}
\op{d}^\bm_{\op{Cac}}(u,v)& \leq \op{d}_{\op{gr}}^m(\rho,u) + \op{d}_{\op{gr}}^m(\rho,v) - 2\Big(\min_{w\in\llbracket u,v\rrbracket} \ell_w - \min \ell -D\Big)\\
&= \ell_u + \ell_v -2 \min_{w\in\llbracket u,v\rrbracket} \ell_w + 2D+2\;.
\end{align*}
This completes the proof. \endproof

\subsection{Random planar maps}

Following \cite{MM} and \cite{MieInvar}, we now discuss Boltzmann
distributions on the space $\m_{r,p}$. We consider a sequence
$\bq=(q_1,q_2,\ldots)$ of non-negative real numbers. We assume that
the sequence $\bq$ has finite support ($q_k=0$ for all sufficiently
large $k$), and is such that $q_k>0$ for some $k\geq 3$. We will then
split our study according to the following two possibilities: 
\begin{enumerate}
\item[(A1)] There exists an odd integer $k$ such that $q_k>0$.
\item[(A2)] The sequence $\bq$ is supported on even integers. 
\end{enumerate}
If $m\in \m_{r,p}$, we define
$$W_\bq(m)= \prod_{f\in F(m)} q_{\op{deg}(f)}$$
where $F(m)$ stands for the set of all faces of $m$ and 
$\op{deg}(f)$ is the degree of the face $f$. In the case when $m=\dagger$, we make
the convention that $q_0=1$ and thus $W_\bq(\dagger)=1$. 

By multiplying the sequence $\bq$ by a suitable positive constant, we
may assume that this sequence is regular critical in the sense of
\cite[Definition 1]{MieInvar} under assumption (A1) or of
\cite[Definition 1]{MaMi} under assumption (A2). We refer the reader
to the Appendix below for details.  In particular, the measure $W_\bq$ is
then finite, and we can define a probability measure $P_\bq$ on
$\m_{r,p}$ by setting
$$P_\bq =Z_\bq^{-1}\,W_\bq,$$
where $Z_\bq=W_\bq(\m_{r,p})$. 

For every integer $n$ such that $W_\bq(\# V(m)=n)>0$,
we consider a random planar map $M_n$ distributed according to the
conditional measure
$$\frac{P_\bq(\cdot \cap \{\# V(m)=n\})}{P_\bq(\# V(m)=n)}.$$
Throughout the remaining part of Section \ref{convcactusmaps},
we restrict our attention to values of $n$ such that $W_\bq(\# V(m)=n)>0$,
so that $M_n$ is well defined. We write $\rho_n$ for the distinguished vertex of $M_n$. 

We now state the main result of this section. In this result, $\bM_n$ stands for the 
graph (pointed at $\rho_n$) 
associated with $M_n$, as explained at the end of subsection \ref{sec:disccac}.

\begin{theorem}
\label{convmapcactus}
There exists a positive constant $B_\bq$ such that
$$B_\bq\,n^{-1/4}\cdot \op{Cac}(\bM_n)
 \build{\la}_{n\to\infty}^{\rm (d)} \op{KAC}$$
 in the Gromov-Hausdorff sense.
 \end{theorem}

The proof of Theorem \ref{convmapcactus} relies on the asymptotic study 
of  the random trees associated with planar maps distributed
under Boltzmann distributions via the BDG bijection. 
The distribution of these random trees was identified in \cite{MaMi} (in the bipartite
case) and in \cite{MieInvar}. 
We set
$$Z^+_\bq= W_\bq (\m_{r,p}^+)\geq 1\quad,\quad Z^0_\bq = W_\bq(\m_{r,p}^0)\,.$$
Note that, under Assumption (A2), $W_\bq$ is supported on bipartite maps and thus
$Z^0_\bq=0$. We also set
$$P_\bq^+=P_\bq(\cdot \mid \m_{r,p}^+)\,,\ P_\bq^-=P_\bq(\cdot \mid \m_{r,p}^-)\,, \ 
P_\bq^0=P_\bq(\cdot \mid \m_{r,p}^0).$$
Note that the definition of $P_\bq^0$ only makes sense under Assumption (A1). 

The next proposition gives the distribution of the tree associated with a 
random planar map distributed according to $P_\bq^+$.
Before stating this proposition,
let us recall that the notion of a four-type Galton-Watson tree is defined 
analogously to the case of a single type. The distribution of such a random
tree is determined by the type of the ancestor, and four offspring
distributions $\nu_i$, $i=1,2,3,4$, which are probability
distributions on $\Z_+^4$; for every $i=1,2,3,4$, $\nu_i$ corresponds to the law of the number of children
(having each of the four possible types) of an individual of type $i$; furthermore, given 
the numbers of children of each type of an individual, these children are ordered
in the tree with the same probability for each possible ordering. See \cite[Section 2.2.1]{MieInvar}
for more details, noting that we consider only the case of  ``uniform ordering'' in the terminology
of \cite{MieInvar}.

\begin{proposition}
\label{BDG-GW}
Suppose that $M^+$ is a random planar map distributed according to
$P_\bq^+$, and let $(\theta, (\L_u)_{u\in\theta_{(1,2)}})$
be the four-type labeled tree associated with $M^+$ via the BDG bijection
between $\TT_+$ and $\m_{r,p}^+$.  Then the
distribution of $(\theta, (\L_u)_{u\in\theta_{(1,2)}})$ is characterized by the following properties:
\begin{enumerate}
\item[\rm(i)] The random tree $\theta$ is a four-type Galton-Watson tree, such that 
the root $\varnothing$ has type $1$ and the offspring distributions $\nu_1,\ldots,\nu_4$
are determined as follows:
\begin{enumerate}
\item[$\bullet$] $\nu_1$ is supported on $\{0\}\times \{0\} \times \Z_+\times\{0\}$, and for every $k\geq 0$,
$$\nu_1(0,0,k,0)= \frac{1}{Z_\bq^+} \Big(1-\frac{1}{Z_\bq^+}\Big)^k.$$
\item[$\bullet$] $\nu_2(0,0,0,1)=1$.
\item[$\bullet$] $\nu_3$ and $\nu_4$ are supported on $\Z_+\times \Z_+\times\{0\}\times\{0\}$, and for every
integers $k,k'\geq 0$,
\begin{align*}
\nu_3(k,k',0,0)&=c_\bq\,(Z^+_\bq)^k(Z^0_\bq)^{k'/2}\,{2k+k'+1\choose k+1} {k+k'\choose k}\,q_{2+2k+k'}\\
 \nu_4(k,k',0,0)&=c'_\bq\,(Z^+_\bq)^k(Z^0_\bq)^{k'/2}\,{2k+k'\choose k} {k+k'\choose k}\,q_{1+2k+k'}
 \end{align*}
 where $c_\bq$ and $c'_\bq$ are the appropriate normalizing constants.
 \end{enumerate}
 \item[\rm(ii)] Conditionally given $\theta$, $(\L_u)_{u\in\theta_{(1,2)}}$ is
 uniformly distributed over all admissible labelings.
 \end{enumerate}
 \end{proposition}
 
 \begin{rek} {\rm The definition of $\nu_4$ does not make sense under Assumption (A2) (because $Z^0_\bq=0$ in that case,
 $\nu_4(k,k',0,0)$ can be nonzero only if $k'=0$, but then $q_{1+2k+k'}=0$). This is however irrelevant since
 under Assumption (A2) the property $Z^0_\bq=0$ entails that $\nu_3$ is supported on $\Z_+\times \{0\}\times\{0\}\times\{0\}$, and thus
 the Galton-Watson tree will have no vertices of type $2$ or $4$. }
 \end{rek}
 
 We refer to \cite[Proposition 3]{MieInvar} for the proof of
 Proposition \ref{BDG-GW} under Assumption (A1) and to
 \cite[Proposition 7]{MaMi} for the case of Assumption (A2). In fact,
 \cite{MieInvar} assumes that $q_k>0$ for some {\em odd} integer
 $k\geq 3$, but the results in that paper do cover the situation
 considered in the present work.
 
 In the next two subsections, we prove Theorem \ref{convmapcactus} under Assumption (A1). 
 The case when Assumption (A2) holds is much easier and will be treated briefly
 in subsection \ref{bipartitecase}.

\subsection{The shuffling operation}
\label{secshuffling}

As already mentioned, we suppose in this section that
 Assumption (A1) holds. 
We consider the random four-type labeled tree $(\theta,(\L_v)_{v\in \theta_{(1,2)}})$ associated with
the planar map $M^+$ via the BDG bijection, as in Proposition \ref{BDG-GW}. 

Our goal is  to investigate the asymptotic 
behavior, when $n$ tends to $\infty$, of the labeled tree 
$(\theta,(\L_v)_{v\in \theta_{(1,2)}})$ conditioned to have $n-1$ vertices of type $1$
(this corresponds to conditioning $M^+$ on the event $\{\# V(M^+)=n\}$).
As already observed in \cite{MieInvar}, a difficulty arises from the 
fact that the label displacements along the tree are not centered, and so 
the results of \cite{MieGW} cannot be applied immediately. To overcome
this difficulty, we will use an idea of \cite{MieInvar}, which consists in introducing
a ``shuffled'' version of the tree $\theta$. In order to explain this, we
need to introduce some notation. 

Let  $\tau$ be a plane tree and $u=(i_1,\ldots,i_p)\in\tau$.
The tree $\tau$ shifted at $u$ is defined by
$$T_u\tau:=\{v=(j_1,\ldots,j_\ell):(i_1,\ldots,i_p,j_1,\ldots,j_\ell)\in\tau\}.$$
Let $k=k_u(\tau)$ be the number of children of $u$ in $\tau$, and, for every
$1\leq i\leq k$, write $u_{(i)}$ for the $i$-th child of $u$. The tree 
$\tau$ reversed at vertex $u$ is the new tree $\tau^*$ characterized by
the properties:
\begin{enumerate}
\item[$\bullet$] Vertices of $\tau^*$ which are not descendants of $u$ are the same as
vertices of $\tau$ which are not descendants of $u$.
\item[$\bullet$] $u\in\tau^*$ and $k_u(\tau^*)=k_u(\tau)=k$.
\item[$\bullet$] For every $1\leq i\leq k$, $T_{u_{(i)}}\tau^* = T_{u_{(k+1-i)}}\tau$.
\end{enumerate}

Our (random) shuffling operation will consist in reversing the tree $\tau$
at every vertex of $\tau$ at an odd generation, with probability $1/2$ for every 
such vertex. We now give a more formal description, which will be needed
in our applications. We keep on considering a (deterministic) plane tree $\tau$.
Let $\mathcal{U}^{\rm o}$ stand for the set of all $u\in\mathcal{U}$ such
that $|u|$ is odd. We consider a collection $(\eps_u)_{u\in\mathcal{U}^{\rm o}}$ of
independent Bernoulli variables with parameter $1/2$. We then define 
a (random) mapping $\sigma:\tau \la \mathcal{U}$ by setting, if 
$u=(i_1,i_2,\ldots,i_p)$, 
$$\sigma(u)=(j_1,j_2,\ldots,j_p)$$
where, for every $1\leq \ell\leq p$,
\begin{enumerate}
\item[$\bullet$] if $\ell$ is odd, $j_\ell=i_\ell$,
\item[$\bullet$] if $\ell$ is even,
$$j_\ell=\left\{\begin{array}{ll}
i_\ell &\hbox{if }\eps_{(i_1,\ldots,i_{\ell-1})} =0\,,\\
k_{(i_1,\ldots,i_{\ell-1})}(\tau)+1 - i_\ell\qquad&\hbox{if }\eps_{(i_1,\ldots,i_{\ell-1})} =1\,.
\end{array}
\right.$$
\end{enumerate}
Then $\tilde \tau =\{\sigma(u): u\in \tau\}$ is a (random) plane tree, called 
the tree derived from $\tau$ by the shuffling operation. If $\tau$ is a
four-type tree, we also view $\tilde \tau$ as a four-type tree by
assigning to the vertex $\sigma(u)$ of $\tilde \tau$ the type of the
vertex $u$ in $\tau$. 

For our purposes it is very important to note that the bijection
$\sigma:\tau \la \tilde\tau$ preserves the genealogical structure, in the sense that 
$u$ is an ancestor of $v$ in $\tau$ if and only if $\sigma(u)$ is an ancestor of $\sigma(v)$
in $\tilde\tau$. Consequently, if $u$ and $v$ are any two vertices of $\tau_{(1,2)}$,
$\llbracket \sigma(u),\sigma(v) \rrbracket$ is the image under $\sigma$
of the set $\llbracket u, v \rrbracket $. 

We can apply this shuffling operation to the random tree $\theta$ (of course
we assume that
the collection $(\eps_u)_{u\in\mathcal{U}^{\rm o}}$ is independent
of $(\theta,(\L_v)_{v\in \theta_{(1,2)}})$). We write $\tilde \theta$ for the four-type tree derived from
$\theta$ by the shuffling operation and we use the same notation $\sigma$
as above for the ``shuffling bijection'' from $\theta$ onto $\tilde\theta$. We assign labels to the vertices 
of $\tilde\theta_{(1,2)}$ by putting for every $u\in\theta_{(1,2)}$,
$$\tilde\L_{\sigma(u)}= \L_u.$$
Note that the random tree $\tilde \theta$ has the same distribution as $\theta$,
and is therefore a four-type Galton-Watson tree as described in 
Proposition \ref{BDG-GW}. On the other hand, the labeled trees
$(\theta,(\L_v)_{v\in \theta_{(1,2)}})$ and 
$(\tilde\theta,(\tilde\L_v)_{v\in\tilde\theta_{(1,2)}})$ have a different
distribution because the admissibility property of labels is not
preserved under the shuffling operation. We can still describe the distribution of the 
labels in the shuffled tree in a simple way. To  this end, write $\op{tp}(u)$ for the type
of a vertex $u$. Then conditionally on $\tilde\theta$, for every vertex $u$
of $\tilde\theta$ such that $|u|$ is odd, if $u_{(1)},\ldots,u_{(k)}$ are the children
of $u$ in lexicographical order, and if $u_{(0)}$ is the parent of $u$, the
vector of label increments
$$(\tilde\L_{u_{(1)}}-\tilde\L_{u_{(0)}},\ldots,\tilde\L_{u_{(k)}}-\tilde\L_{u_{(0)}})$$ is with probability $1/2$ uniformly distributed
over the set
$$\mathbb{A}:=\{(i_1,\ldots,i_k)\in \Z^k: i_{j+1}\geq i_j -\mathbf{1}_{\{\op{tp}(u_{(j+1)})=1\}}\;,\hbox{ for all }0\leq j\leq k\},$$
and with probability $1/2$ uniformly distributed
over the set
$$\mathbb{A}':=\{(i_1,\ldots,i_k)\in \Z^k: i_{j}\geq i_{j+1}-\mathbf{1}_{\{\op{tp}(u_{(j)})=1\}}\;,\hbox{ for all }0\leq j\leq k\}.$$
In the definition of both $\mathbb{A}$ and $\mathbb{A}'$ we make the convention that $i_0=i_{k+1}=0$ and $u_{(k+1)}=u_{(0)}$.
Furthermore the vectors of label increments are independent (still conditionally on $\tilde\theta$) when $u$ varies over 
vertices of $\tilde \theta$ at odd generations.

The preceding description of the distribution of labels in the shuffled tree is easy to
establish. Note that the set $\mathbb{A}$ corresponds to the admissibility property of
labels, whereas $\mathbb{A}'$ corresponds to a ``reversed'' version of this 
property. 

For every $u\in \tilde\theta_{(1,2)}$, set
$$\tilde\L_u' = \tilde \L_u -\frac{1}{2}\mathbf{1}_{\{\op{tp}(u)=2\}}.$$
If we replace $\tilde \L_u$ by $\tilde\L_u'$, then the vectors of label increments in
$\tilde\theta$ become centered. This follows from elementary arguments: See \cite[Lemma 2]{MieInvar}
for a detailed proof. As in \cite{MieInvar} or in \cite{MW}, the fact that the label increments are centered allows us to use the asymptotic results of
\cite{MieGW}, noting that these results will apply to $\tilde \L_u$ as well as to $\tilde\L_u'$ 
since the additional term $\frac{1}{2}\mathbf{1}_{\{\op{tp}(u)=2\}}$ obviously plays no role in the scaling limit. 
Before we state the relevant result, we need to introduce some notation.

For $n\geq 2$, let $(\tilde\theta^n,(\tilde\L^n_v)_{v\in \tilde\theta^n_{(1,2)}})$ be distributed as the
labeled tree $(\tilde\theta,(\tilde\L_v)_{v\in \tilde\theta_{(1,2)}})$ conditioned on the event
$\{\#\tilde\theta_{(1)}=n-1\}$ (recall that we restrict our attention to values of $n$
such that the latter event has positive probability).
Let $p_n=\#\tilde\theta^n-1$ and let $u^n_0=\varnothing, u^n_1,\ldots,u^n_{p_n}=\varnothing$ be the modified contour sequence
of $\tilde\theta_n$.
The contour process $C^n=(C^n_i)_{0\leq i\leq p_n}$ is defined by
$$C^n_i = |u^n_i|$$
and the label process $V^n=(V^n_i)_{0\leq i\leq p_n}$ by
$$V^n_i = \tilde\L^n_{u^n_i}\;.$$
We extend the definition of both processes $C^n$ and $V^n$ to the real interval $[0,p_n]$
by linear interpolation. 

Recall the notation $(\be,Z)$ from Section \ref{Brcactus}.

\begin{proposition}
\label{convcoding}
There exist two positive constants $A_\bq$ and $B_\bq$ such that
\begin{equation}
\label{keyconv}
\Big(A_\bq\,\frac{C^n(p_n s)}{n^{1/2}}, B_\bq\,\frac{V^n(p_ns)}{n^{1/4}}\Big)_{0\leq s\leq 1}
\build{\la}_{n\to\infty}^{\rm (d)} (\be_s, Z_s)_{0\leq s\leq 1}
\end{equation}
in the sense of weak convergence of the distributions on the space $C([0,1],\R^2)$. 
\end{proposition}

This follows from the more general results proved in \cite{MieGW} for spatial
mutitype Galton-Watson trees. One should note that the results of
\cite{MieGW} are given for variants of the contour process and the label process
(in particular the contour process is replaced by the so-called height process
of the tree). However simple arguments show that the convergence in the proposition
can be deduced from the ones in \cite{MieGW}: See in particular Section 1.6 of \cite{tree}
for a detailed explanation of why convergence results for the height process 
imply similar results for the contour process. Proposition
\ref{convcoding} is also equivalent to Theorem 3.1 in \cite{MW},
where the contour and label processes are defined in a slightly different way.

\subsection{Proof of Theorem \ref{convmapcactus} under Assumption (A1)}
\label{proofMainT}

We keep assuming that Assumption (A1) holds. 
Let $M_n^+$ be distributed according to the probability measure $P_\bq^+(\cdot\mid \# V(m)=n)$,
or equivalently as $M^+$ conditionally on the event $\{\# V(M^+)=n\}$. As above, $\rho_n$
stands for the distinguished point of $M_n^+$, and we
will write $\bM_n^+$ for the pointed graph associated with $M_n^+$. Let $(\theta^n,(\L^n_v)_{v\in \theta^n_{(1,2)}})$
be the random labeled tree associated with $M_n^+$ via the BDG bijection
between $\TT_+$ and $\mathcal{M}_{r,p}^+$. Notice that 
$(\theta^n,(\L^n_v)_{v\in \theta^n_{(1,2)}})$ has the same distribution as 
$(\theta,(\L_v)_{v\in \theta_{(1,2)}})$ conditional on $\{\#\theta_{(1)}=n-1\}$. 

\smallskip
We write $(\tilde\theta^n,(\tilde\L^n_v)_{v\in \tilde\theta^n_{(1,2)}})$
for the tree derived from $(\theta^n,(\L^n_v)_{v\in\theta^n_{(1,2)}})$ by the shuffling operation, and
$\sigma_n$ for the shuffling bijection from $\theta^n$ onto $\tilde\theta^n$. 
The notation $(\tilde\theta^n,(\tilde\L^n_v)_{v\in\tilde\theta^n_{(1,2)}})$ is consistent with
the end of the preceding subsection, since conditioning  the tree on having $n-1$ vertices of type $1$
clearly commutes with the shuffling operation.

\smallskip
As previously,
$u^n_0=\varnothing,u^n_1,\ldots,u^n_{p_n}$ denotes the modified contour sequence of $\tilde\theta^n$.
For every $j\in\{0,1,\ldots,p_n\}$, we put $v^n_j=\sigma_n^{-1}(u^n_j)$. Recall that 
by construction the type of $u^n_j$ (in $\tilde \theta^n$) coincides with the type
of $v^n_j$ (in $\theta^n$).

  Using the Skorokhod representation theorem, we may assume that the
 convergence (\ref{keyconv}) holds almost surely. We will then prove that the
 convergence  
\begin{equation}
\label{convmapcactus1}
B_\bq\,n^{-1/4}\cdot \op{Cac}(\bM^+_n)
 \build{\la}_{n\to\infty}^{} \op{KAC}
 \end{equation}
also holds almost surely, in the Gromov-Hausdorff sense. 
 
 We first define a correspondence $\mathcal{R}_n^0$ between $\T_\be$ and $V(M^+_n)$
 by declaring that $(a_*,\rho_n)$ belongs to $\mathcal{R}_n^0$, and, for every $s\in[0,1]$:
 \begin{enumerate}
 \item[$\bullet$] if $v^n_{[p_ns]}$ is of type $1$, $(p_\be(s),v^n_{[p_ns]})$ belongs to $\mathcal{R}_n^0$;
 \item[$\bullet$] if $v^n_{[p_ns]}$ is of type $2$, then if $w$ is any of the two (possibly equal) vertices of type $1$ associated with
 $v^n_{[p_ns]}$, $(p_\be(s),w)$ belongs to $\mathcal{R}_n^0$.
 \end{enumerate}
 We then write $\mathcal{R}_n$ for the induced correspondence between the quotient spaces 
 $\op{KAC}=\T_\be\,/\!\asymp$ and $ \op{Cac}(\bM_n^+)$. A pair $(x,\alpha)\in\op{KAC}\times \op{Cac}(\bM_n^+)$ belongs to $\mathcal{R}_n$
 if and only if there exists a representative $a$ of $x$ in $\T_\be$ and a representative 
 $u$ of $\alpha$ in $V(M_n^+)$ such that $(a,u)\in\mathcal{R}_n^0$. 

Thanks to (\ref{GHcorres}), the convergence (\ref{convmapcactus1}) will be proved if we can verify that the distortion of $\mathcal{R}_n$,
when $\op{KAC}$ is equipped with the distance $\op{d}_{\op{KAC}}$ and $ \op{Cac}(\bM_n^+)$ is equipped 
with $ B_\bq\,n^{-1/4}\,\op{d}^{\bM_n^+}_{\op{Cac}}$, 
tends to $0$ as $n\to\infty$, almost surely. To this end, it is enough to verify that
\begin{equation}
\label{distorttech1}
\lim_{n\to \infty} \sup_{0\leq s\leq 1} \Big| \op{d}_{\op{KAC}}(a_*,p_\be(s)) - B_\bq\,n^{-1/4}\,\op{d}^{\bM_n^+}_{\op{Cac}}(\rho_n,\widehat v^n_{[p_ns]})\Big|
=0\;,\quad\hbox{a.s.}
\end{equation}
and 
\begin{equation}
\label{distorttech2}
\lim_{n\to\infty}\sup_{s,t\in[0,1]} \Big| \op{d}_{\op{KAC}}(p_\be(s),p_\be(t)) - B_\bq 
\,n^{-1/4}\,\op{d}^{\bM_n^+}_{\op{Cac}}(\widehat v^n_{[p_ns]},\widehat v^n_{[p_nt]})\Big| = 0 \;,\quad\hbox{a.s.}
\end{equation}
In both (\ref{distorttech1}) and (\ref{distorttech2}), $\widehat v^n_{[p_ns]}=v^n_{[p_ns]}$ if $v^n_{[p_ns]}$
is of type $1$, whereas, if $v^n_{[p_ns]}$ is of type $2$, $\widehat v^n_{[p_ns]}$ stands for one of the vertices
of type $1$ associated with $v^n_{[p_ns]}$ (obviously the validity of (\ref{distorttech1}) and (\ref{distorttech2})
does not depend on the choice of this vertex). 

The proof of (\ref{distorttech1}) is easy. Note that
$$\op{d}_{\op{KAC}}(a_*,p_\be(s)) = Z_{p_\be(s)}-Z_{a_*}= Z_s - \underline Z$$
and, by (\ref{distance-label}),
$$\op{d}^{\bM_n^+}_{\op{Cac}}(\rho_n,\widehat v^n_{[p_ns]})
= \op{d}_{\op{gr}}^{M_n^+}(\rho_n,\widehat v^n_{[p_ns]})= \L^n_{\widehat v^n_{[p_ns]}}-\min\L^n +1$$
so that
$$| \op{d}^{\bM_n^+}_{\op{Cac}}(\rho_n,\widehat v^n_{[p_ns]}) - (\L^n_{v^n_{[p_ns]}}-\min\L^n)|\leq 1.$$
Since $\L^n_{v^n_{[p_ns]}}-\min\L^n= \tilde\L^n_{u^n_{[p_ns]}}-\min\tilde\L^n=V^n_{[p_ns]} - \min V^n$, our claim (\ref{distorttech1})
follows from the (almost sure) convergence (\ref{keyconv}).

It remains to establish (\ref{distorttech2}). It suffices to prove that
almost surely, for every choice of the sequences $(s_n)$ and $(t_n)$ in $[0,1]$,
we have
$$\lim_{n\to\infty}  \Big| \op{d}_{\op{KAC}}(p_\be(s_n),p_\be(t_n)) - B_\bq 
\,n^{-1/4}\,\op{d}^{\bM_n^+}_{\op{Cac}}(\widehat v^n_{[p_ns_n]},\widehat v^n_{[p_nt_n]})\Big| = 0.$$
We will prove that the preceding convergence holds 
for all choices of the sequences $(s_n)$ and $(t_n)$, on the set of full probability measure
where the convergence (\ref{keyconv}) holds. From now on we
argue on the latter set.

By a compactness argument, we may assume that the sequences $(s_n)$ and $(t_n)$
converge to $s$ and $t$ respectively
as $n\to\infty$. The proof then reduces to checking that
$$\lim_{n\to\infty} B_\bq 
\,n^{-1/4}\,\op{d}^{\bM_n^+}_{\op{Cac}}(\widehat v^n_{[p_ns_n]},\widehat v^n_{[p_nt_n]})
= \op{d}_{\op{KAC}}(p_\be(s),p_\be(t))= Z_s+Z_t - 2\,\min_{c\in\llbracket p_\be(s),p_\be(t)\rrbracket} Z_c\,.$$
From Corollary \ref{bounddiscretecactus} (and the fact that the sequence $\bq$ is finitely supported), this
will follow if we can verify that
$$\lim_{n\to\infty} B_\bq 
\,n^{-1/4}\, \Big(\L^n_{\widehat v^n_{[p_ns_n]}} + \L^n_{\widehat v^n_{[p_nt_n]}} -2 \min_{w\in \llbracket \widehat v^n_{[p_ns_n]},
\widehat v^n_{[p_nt_n]}\rrbracket} \L^n_w\Big) = Z_s+Z_t - 2\,\min_{c\in\llbracket p_\be(s),p_\be(t)\rrbracket} Z_c\,.$$
Observe that 
$$|\L^n_{\widehat v^n_{[p_ns_n]}} - \L^n_{v^n_{[p_ns_n]}}| \leq 1$$
and $\L^n_{v^n_{[p_ns_n]}} = \tilde \L^n_{u^n_{[p_ns_n]}}$.
From the convergence (\ref{keyconv}), we have
$$\lim_{n\to\infty} B_\bq 
\,n^{-1/4}\,\L^n_{\widehat v^n_{[p_ns_n]}}= \lim_{n\to\infty} B_\bq 
\,n^{-1/4}\,\tilde\L^n_{u^n_{[p_ns_n]}}= \lim_{n\to\infty} B_\bq 
\,n^{-1/4}\,V^n_{[p_ns_n]} = Z_s$$
and similarly if the sequence $(s_n)$ is replaced by $(t_n)$. Finally, we need to verify that
\begin{equation}
\label{mapcactustech}
\lim_{n\to\infty} \Big(B_\bq\,n^{-1/4}\,\min_{w\in \llbracket \widehat v^n_{[p_ns_n]},
\widehat v^n_{[p_nt_n]}\rrbracket} \L^n_w\Big) = 
\min_{c\in\llbracket p_\be(s),p_\be(t)\rrbracket} Z_c\,.
\end{equation}
In proving (\ref{mapcactustech}), we may replace $\widehat v^n_{[p_ns_n]}$ 
and $\widehat v^n_{[p_nt_n]}$ by $v^n_{[p_ns_n]}$, 
and $v^n_{[p_nt_n]}$ respectively. The point is that if $u$ is a vertex of $\theta^n$ of type $2$
and $v$ is an associated vertex of type $1$, our definitions imply that $\min_{w\in\llbracket u,v\rrbracket} \L^n_w = \L^n_v$.
Without loss of generality we can also assume that $s\leq t$. 

Since $\llbracket u^n_{[p_ns_n]},
u^n_{[p_nt_n]}\rrbracket$ is the image under $\sigma_n$ of
$\llbracket v^n_{[p_ns_n]},
v^n_{[p_nt_n]}\rrbracket$, (\ref{mapcactustech})
will hold if we can prove that
\begin{equation}
\label{mapcactustechh}
\lim_{n\to\infty} \Big(B_\bq\,n^{-1/4}\,\min_{w\in \llbracket u^n_{[p_ns_n]},
u^n_{[p_nt_n]}\rrbracket} \tilde\L^n_w\Big) = 
\min_{c\in\llbracket p_\be(s),p_\be(t)\rrbracket} Z_c\,.
\end{equation}
Let us first prove the upper bound
\begin{equation}
\label{mapcactustech1}
\limsup_{n\to\infty} \Big(B_\bq\,n^{-1/4}\,\min_{w\in \llbracket u^n_{[p_ns_n]},
u^n_{[p_nt_n]}\rrbracket} \tilde\L^n_w\Big) \leq
\min_{c\in\llbracket p_\be(s),p_\be(t)\rrbracket} Z_c\,.
\end{equation}
Let us pick $c\in \llbracket p_\be(s),p_\be(t)\rrbracket$. We may assume that
$c\not = p_\be(s)$ and $c\not =p_\be(t)$ (otherwise the desired lower bound immediately
follows from the convergence (\ref{keyconv})). Then, we can find 
$r\in(s,t)$ such that $c=p_\be(r)$ and either
$$\be_u> \be_r\;,\quad\hbox{for every }u\in[s,r)$$
or 
$$\be_u > \be_r\;,\quad\hbox{for every }u\in (r,t].$$
Consider only the first case, since the second one can be treated in a similar manner.
The convergence of the rescaled contour processes 
then guarantees that we can find a sequence $(k_n)$ of positive integers
such that $k_n/p_n \la r$ as $n\to\infty$, and 
$$C^n_k > C^n_{k_n}\;,\quad\hbox{for every }k\in\{[p_ns_n],[p_ns_n]+1,\ldots,k_n-1\}$$
for all sufficiently large $n$. The latter property, and the construction of the contour
sequence of the tree $\theta^n$, ensure that $u^n_{k_n}\in \llbracket u^n_{[p_ns_n]},u^n_{[p_nt_n]}\rrbracket$, for all
sufficiently large $n$.
However, by the convergence of the rescaled label processes, we have
$$\lim_{n\to\infty} B_\bq\,n^{-1/4}\, \tilde\L^n_{u^n_{k_n}}= Z_r = Z_c.$$
Consequently,
$$\limsup_{n\to\infty} \Big(B_\bq\,n^{-1/4}\,\min_{w\in \llbracket u^n_{[p_ns_n]},
u^n_{[p_nt_n]}\rrbracket} \tilde\L^n_w\Big) \leq Z_c$$
and since this holds for every choice of $c$
the upper bound (\ref{mapcactustech1}) follows.

Let us turn to the lower bound
\begin{equation}
\label{mapcactustech2}
\liminf_{n\to\infty} \Big(B_\bq\,n^{-1/4}\,\min_{w\in \llbracket u^n_{[p_ns_n]},
u^n_{[p_nt_n]}\rrbracket} \tilde\L^n_w\Big) \geq
\min_{c\in\llbracket p_\be(s),p_\be(t)\rrbracket} Z_c\,.
\end{equation}
For every $n$, let $w_n\in \llbracket u^n_{[p_ns_n]},
u^n_{[p_nt_n]}\rrbracket$ be such that
$$\min_{w\in \llbracket u^n_{[p_ns_n]},
u^n_{[p_nt_n]}\rrbracket} \tilde\L^n_w=\tilde\L^n_{w_n}.$$
We can write $w_n=u^n_{j_n}$ where $j_n\in \{[p_ns_n],[p_ns_n]+1,\ldots,[p_nt_n]\}$
is such that 
\begin{equation}
\label{mapcactustech3}
C^n_{j_n}= \min_{[p_ns_n]\leq j\leq j_n} C^n_j\;,
\end{equation}
or
\begin{equation}
\label{mapcactustech4}
C^n_{j_n}= \min_{j_n\leq j\leq [p_nt_n]} C^n_j\;.
\end{equation}
We need to verify that
$$\liminf_{n\to\infty} B_\bq\,n^{-1/4}\,\tilde\L^n_{w_n}\geq 
\min_{c\in\llbracket p_\be(s),p_\be(t)\rrbracket} Z_c\,.$$
We argue by contradiction and suppose that there exist $\varepsilon>0$
and a subsequence $(n_k)$ such that, for every $k$,
$$B_\bq\,n_k^{-1/4}\,\tilde\L^{n_k}_{w_{n_k}}\leq 
\min_{c\in\llbracket p_\be(s),p_\be(t)\rrbracket} Z_c -\varepsilon.$$
By extracting another subsequence if necessary, we may assume 
furthermore that $j_{n_k}/p_{n_k} \la r\in[s,t]$ as $k\to\infty$, and that
(\ref{mapcactustech3}) holds with $n=n_k$ for every $k$ (the case when 
(\ref{mapcactustech4}) holds instead of (\ref{mapcactustech3}) is treated in a similar
manner). Then, from the convergence of rescaled contour processes,
we have 
$$\be_r=\min_{s\leq u\leq r} \be_r\,,$$
which implies that $p_\be(r)\in\llbracket p_\be(s),p_\be(t)\rrbracket$. 
Furthermore, from the convergence of rescaled label processes,
$$Z_{p_\be(r)}=Z_r =\lim_{k\to\infty} B_\bq\,n_k^{-1/4}\,\tilde\L^{n_k}_{w_{n_k}}
\leq \min_{c\in\llbracket p_\be(s),p_\be(t)\rrbracket} Z_c -\varepsilon.$$
This contradiction completes the proof of (\ref{mapcactustech2}) and
of the convergence (\ref{convmapcactus1}).

In order to complete the proof of Theorem \ref{convmapcactus} under Assumption (A1),
it suffices to verify that the convergence (\ref{convmapcactus1}) also holds (in distribution)
if $M_n^+$ is replaced by a random planar map $M_n^-$ distributed
according to $P^-_\bq(\cdot \mid \# V(m)=n)$, or by a random planar map $M_n^0$ distributed
according to $P^0_\bq(\cdot \mid \# V(m)=n)$. The first case is trivial since $M_n^-$
can be obtained from $M_n^+$ simply by reversing the orientation of the root edge. The case
of $M_n^0$ is treated by a similar method as the one we used for $M_n^+$. We first need
an analogue of Proposition \ref{BDG-GW}, which is provided by the last statement of 
Proposition 3 in \cite{MieInvar}. In this analogue, the random labeled tree associated
with a planar map distributed according to $P_\bq^0$ is described as the concatenation 
(at the root vertex) of two independent labeled Galton-Watson trees whose root is of type $2$, with the
same offspring distributions as in Proposition \ref{BDG-GW}. The results of \cite{MieGW}
can be used to verify that Proposition \ref{convcoding} still holds with the
same constants $A_\bq$ and $B_\bq$, and the remaining part of the argument
goes through without change. This completes the proof of
Theorem \ref{convmapcactus} under Assumption (A1). 

\subsection{The bipartite case}
\label{bipartitecase}

In this section, we briefly discuss the proof of Theorem \ref{convmapcactus} under Assumption (A2). 
In that case, since $W_\bq(\mathcal{M}^0_{r,p})=0$, it is obviously enough to prove the convergence 
of Theorem \ref{convmapcactus} with $M_n$ replaced by $M_n^+$. 
The proof becomes
much simpler because we do not need the shuffling operation. As previously, we
introduce the labeled tree $(\theta^n,(\L^n_v)_{v\in \theta^n_{(1,2)}})$
associated with $M_n^+$ via the BDG bijection, but we  now define
$u^n_0=\varnothing, u^n_1,\ldots,u^n_{p_n}=\varnothing$ as the modified contour sequence
of $\theta_n$ (instead of $\tilde \theta_n$). We then define the contour process $C^n_i=|u^n_i|$
and the label process $V^n_i=\L^n_{u^n_i}$, for $0\leq i\leq p_n$. Proposition 3.7 then holds in
exactly the same form, as a consequence of the results of \cite{MaMi}. The reason why we do not need
the shuffling operation is the fact that the label increments of $(\theta^n,(\L^n_v)_{v\in \theta^n_{(1,2)}})$
are centered in the bipartite case.

Once the convergence (\ref{keyconv}) is known to hold, it suffices to repeat all steps of the 
proof in subsection \ref{proofMainT}, replacing $\tilde\theta^n$ by $\theta^n$ and $v^n_i$ by $u^n_i$
wherever this is needed. We leave the details to the reader.

\section{The dimension of the Brownian cactus}
\label{hauscalc}

In this section, we compute the Hausdorff dimension of the 
Brownian cactus $\op{KAC}$. We write $\mathfrak{p}:\T_\be \la\op{KAC}= \T_\be\,/\!\asymp$ for the canonical
projection. The uniform measure $\mu$ on $\op{KAC}$ is the image of the mass measure $\op{Vol}$ on the CRT (see
Section \ref{Brcactus}) under $\mathfrak{p}$. 
For every $x$ in $\op{KAC}$ and every $\delta\geq 0$, we denote  the closed ball of center $x$ and radius $\delta$ in $\op{KAC}$ 
by $B_{\op{KAC}}(x,\delta)$. The following theorem gives information about the $\mu$-measure of these balls around a typical point of $\op{KAC}$. 

\begin{proposition} \label{haus} 
{\rm(i)} We have 
$$\displaystyle\E{ \int \mu(dx)\,\mu \big(B_{\op{KAC}}(x,\delta)\big)} = \frac{2^{5/4}\,\Gamma(1/4)}{3\sqrt{\pi}}\;\delta^3 +o(\delta^3),$$
as $\delta\to 0$.

\noindent{\rm(ii)} For every $\varepsilon >0$, 
$$\displaystyle\limsup_{\delta \to 0} \frac{\mu\big(B_{\op{KAC}}(x,\delta)\big)}{\delta^{4-\varepsilon}} = 0\;,\ \mu(dx)\ \hbox{a.e.}, \ \hbox{a.s.}$$
\end{proposition}

\begin{rek} {\rm Let $U$ be uniformly distributed over $[0,1]$, so that $p_\be(U)$ is distributed according to $\op{Vol}$
and $X=\mathfrak{p}\circ p_\be(U)$ is distributed according to $\mu$. Assertion (i) of the theorem says that the mean volume of
the ball $B_{\op{KAC}}(X,\delta)$ is of order $\delta^3$, whereas
assertion (ii) shows that almost surely the volume of this ball 
will be bounded above by $\delta^{4-\varepsilon}$ when $\delta$ is small. This 
difference between the mean and the almost sure behavior is 
specific to the Brownian cactus. In the case of the Brownian map, results
from Section 6 of \cite{Geo} show that $\delta^4$ is the correct order both for 
the mean and the almost sure behavior of the volume of
a typical ball of radius $\delta$. 

In relation with this, we see that in contrast with the CRT or the
Brownian map, the Brownian cactus is not invariant
under re-rooting according to the ``uniform'' measure $\mu$. This
means that $\op{KAC}$ re-rooted at $X$ does not have the same
distribution as $\op{KAC}$. Indeed, since
$\op{d}^{\mathbf{E}}_{\op{Kac}}(\rho,x)=d(\rho,x)$ for every pointed
geodesic space $\mathbf{E}=(E,d,\rho)$, the previous considerations, 
and Proposition \ref{identKactus},
entail that $\mu(B_{\op{KAC}}(\rho,\delta))$ is of order $\delta^4$
both in the mean and in the a.s.\ sense. }
\end{rek}

\proof (i) Fix $\delta>0$. Let $U$ and $U'$ be two independent random variables that are uniformly
distributed over $[0,1]$ and independent of $(\be,Z)$. By the very definition of $\mu$, we have 
$$\E{ \int \mu(dx)\,\mu \big(B_{\op{KAC}}(x,\delta)\big)} = \P{\op{d}_{\op{KAC}}(p_{\be}(U),p_{\be}(U')) \leq \delta}.$$ 
The value of $\op{d}_{\op{KAC}}(p_{\be}(U),p_{\be}(U'))$ is determined by the labels $Z_{a}$ for $a \in \llbracket p_{\be}(U),p_{\be}(U')\rrbracket$. 
Write $(g_{U,U'}(t),0\leq t\leq \op{d}_\be(U,U'))$ for the geodesic path from $p_\be(U)$ to $p_\be(U')$ in the tree $\T_\be$
(so that $\llbracket p_\be(U),p_\be(U')\rrbracket$ is the range of $g_{U,U'}$). 
Then, conditionally on the triplet $(\be,U,U')$ the process 
$$ \Big(Z_{g_{U,U'}(t)}-Z_{p_{\be}(U)}\Big)_{0\leq t \leq d_{\be}(U,U')},$$ is a standard linear Brownian motion. Hence 
if $(B_{t})_{t\geq0}$ is a linear Brownian motion independent of  $(\be,U,U')$, we have
\begin{eqnarray*} \P{\op{d}_{\op{KAC}}(p_\be(U),p_\be(U') )\leq \delta} &=& \P{B_{L}-2 \min_{0\leq s\leq L} B_{s} \leq \delta}, \end{eqnarray*}
where $L= \op{d}_\be (U,U')$. 
Pitman's theorem \cite[Theorem VI.3.5]{RY} implies that, for every fixed $l\geq 0$, $B_{l}-2\min_{0\leq s\leq l} B_{s} $ has the same
distribution as $B^{(3)}_{l}$, where $(B^{(3)}_{t})_{t \geq 0}$ denotes  a three-dimensional Bessel process started from $0$. 
From the invariance under uniform re-rooting of the distribution of the CRT (see for example \cite{LGW}), the variable $\op{d}_{\be}(U,U')$ has the same distribution as $\op{d}_{\be}(0,U)=\be_U$, which has density $4l\,e^{-2l^2}$. Consequently, we can explicitly compute 
\begin{eqnarray*}
\P{\op{d}_{\op{KAC}}(U,U') \leq \delta} &=& 4\displaystyle\int_{0}^\infty dl\;le^{-2l^2} \P{B^{(3)}_{l} \leq \delta},\\
&=&4 \displaystyle\int_{0}^\infty dl\;l e^{-2l^2} \int_{\R^3} dz\;(2\pi l)^{-3/2}\,e^{-|z|^2/2l}\,\mathbf{1}_{\{|z|\leq \delta\}},\\
 &=&\displaystyle4\sqrt{\frac{2}{\pi}}\int_{0}^\infty dl\;l^{-1/2} e^{-2l^2}\int_{0}^\delta du\;u^2\, e^{-u^2/2l},\\
 &=&\displaystyle  4\sqrt{\frac{2}{\pi}} \int_{0}^\delta du \;u^2\int_{0}^{\infty} dl \;l^{-1/2} \exp\left(-2l^2-(u^2/2l)\right).
  \end{eqnarray*}
  The desired result follows since
  $$\lim_{u\to 0} \int_{0}^{\infty} dl \;l^{-1/2} \exp\left(-2l^2-(u^2/2l)\right) = 
  \int_{0}^{\infty} dl \;l^{-1/2} \exp\left(-2l^2\right)= 2^{-5/4}\,\Gamma(1/4). $$

(ii) Let us fix $r\in]0,1[$. For every $u\in[0,\be_r]$, set
\begin{align*}
&G_\be(r,u) = \max\{ s \in [0,r]:\be_s= \be_r-u\},\\
&D_\be(r,u) = \min\{s \in [r,1]:\be_s=\be_r-u\}.
\end{align*}
Then $p_\be(G_\be(r,u))=p_\be(D_\be(r,u))$ is a point of $\llbracket p_\be(0),p_\be(r)\rrbracket$, and more precisely
the path $u\la p_\be(G_\be(r,u))$, $0\leq u\leq \be_r$ is the geodesic from $p_\be(r)$ to $p_\be(0)$
in the tree $\T_\be$. As a consequence, conditionally on $\be$, the process
$$M^{(r)}_u:=Z_r - \min\{Z_v: v\in \llbracket p_\be(G_\be(r,u)),p_\be(r)\rrbracket\}\;,\quad 0\leq u\leq \be_r$$
has the same distribution as 
$$-\min_{0\leq v\leq u} B_v\;,\quad 0\leq u\leq \be_r$$
where $B$ is as above. By classical results
(see e.g. Theorem 6.2 in \cite{JT}), we have, for every $\varepsilon\in]0,1/2[$,
\begin{equation}
\label{HDtech}
\lim_{u \to 0} u^{-1/2-\varepsilon} M^{(r)}_u= \infty\;,\quad\hbox{a.s.} 
\end{equation}
On the other hand, if $t\in[0,1]\backslash]G_\be(r,u),D_\be(r,u)[$, we have $\min_{t\wedge r\leq s\leq t\vee r} \be_s \leq \be_r-u$,
which implies that the segment $\llbracket p_\be(t),p_\be(r)\rrbracket$ contains 
$\llbracket p_\be(G_\be(r,u)), p_\be(r)\rrbracket$, and therefore
$$\op{d}_{\op{KAC}}(p_\be(t),p_\be(r))\geq M^{(r)}_u.$$
Using (\ref{HDtech}), it follows that, for every fixed $\varepsilon\in]0,1/2[$, we have a.s. for all $u>0$ small enough
$$B_{\op{KAC}}( p_\be(r),u^{1/2 +\varepsilon}) \subset  \Big(\op{KAC}\setminus\, \mathfrak{p}\circ p_\be
\left([0,G_\be(r,u)]\cup [D_\be(r,u),1]\right)
\Big),$$
and in particular
$$\mu(B_{\op{KAC}}( p_\be(r),u^{1/2 +\varepsilon})) \leq D_\be(r,u)-G_\be(r,u).$$
However, the same standard results about Brownian motion that we already used to derive (\ref{HDtech}) imply that
$$ \lim _{u \to 0}u^{-2+\varepsilon} (D_\be(r,u)-G_\be(r,u)) = 0\;,\quad\hbox{a.s.}$$
We conclude that, for every $\varepsilon\in]0,1/2[$,
$$\lim_{u\to 0} u^{-2+\varepsilon}\mu(B_{\op{KAC}}( p_\be(r),u^{1/2 +\varepsilon}))=0\;,\quad\hbox{a.s.}$$
and property (ii) follows, in fact in a slightly stronger form than stated in the theorem. \endproof

\begin{corollary} Almost surely, the Hausdorff dimension of $\op{KAC}$ is $4$.
\end{corollary}
\proof Classical density theorems for Hausdorff measures show that the existence of a non-trivial measure $\mu$
satisfying the property stated in part (ii) of Proposition \ref{haus} implies the lower bound
$\op{dim}(\op{KAC})\geq 4$. To get the corresponding upper bound, we first note that the
mapping $[0,1]\ni t\la Z_t$ is a.s. H\"older continuous with exponent $1/4 - \varepsilon$, for any
$\varepsilon\in]0,1/4[$. Observing that $\llbracket p_{\be}(t),p_{\be}(t')\rrbracket \subset p_{\be}([t\wedge t',t\vee t'])$,
for every $t,t'\in[0,1]$, it readily follows that the composition $\mathfrak{p}\circ p_\be$ defined on $[0,1]$ and with values in 
$\op{KAC}$, is a.s. H\"older continuous with exponent $1/4 - \varepsilon$, for any
$\varepsilon\in]0,1/4[$. Hence, the Hausdorff dimension of $\op{KAC}$, which is
the range of $\mathfrak{p}\circ p_\be$, must be bounded above by $4$. \endproof

\section{Separating cycles}
\label{cycle}

In this section, we study the existence and properties of a cycle with minimal length 
separating two points of the Brownian map, under the condition that this cycle contains 
a third point. This is really a problem about the Brownian map, but the cactus distance 
plays an important role in the statement. Our results in this section are related to
the work of Bouttier and Guitter \cite{BG} for large random quadrangulations of the plane.

We consider the Brownian map as the random pointed compact metric space 
$(m_\infty, D,\rho_*)$ that appears in the convergence (\ref{conv-to-Bmap}) for a suitable
choice of the sequence $(n_k)$. Recall that the metric $D$ may depend
on the choice of the sequence, but the subsequent results will hold for any of
the possible limiting metrics. We set $\bp=\Pi\circ p_\be$, which corresponds to the
canonical projection from $[0,1]$ onto $m_\infty$. If $U$ is uniformly distributed
over $[0,1]$, the point $\bp(U)$ is distributed according to the volume measure
$\lambda$ on $m_\infty$.

A loop in $m_\infty$ is a continuous path $\gamma:[0,T]\la m_\infty$, where $T>0$, such that
$\gamma(0)=\gamma(T)$. If $x$ and $y$ are two distinct points of $m_\infty$, we say that the loop $\gamma$ separates the points $x$
and $y$ if $x$ and $y$ lie in distinct connected components of $m_\infty\backslash \{\gamma(t):0\leq t\leq T\}$.
It is known \cite{LGP} that $(m_\infty, D)$ is homeomorphic to the $2$-sphere, so that 
separating loops do exist. We denote by $S(x,y,\rho_*)$ the set of all loops
$\gamma$ such that $\gamma(0)=\rho_*$ and $\gamma$ separates $x$ and $y$.
Recall from subsection \ref{sec:contiK} the definition of the length of a curve
in a metric space.

\begin{theorem}
\label{separating-cycle}
Let $U_1$ and $U_2$ be independent and uniformly distributed over 
$[0,1]$. Then almost surely there exists  a unique loop $\gamma_*\in S(\bp(U_1),\bp(U_2),\rho_*)$
with minimal length, up to reparametrization and time-reversal. This loop is obtained as the concatenation of the two distinct geodesic 
paths from $\Pi(\beta)$ to $\rho_*$, where $\beta$ is the a.s. unique point of 
$\llbracket p_\be(U_1),p_\be(U_2)\rrbracket$ such that
$$Z_\beta=\min_{a\in \llbracket p_\be(U_1),p_\be(U_2)\rrbracket} Z_a.$$
In particular, the length of $\gamma_*$ is
$$L(\gamma_*)= 2 D(\rho_*,\Pi(\beta))=D(\rho_*,\bp(U_1)) + D(\rho_*,\bp(U_2)) - 2\,\op{d}_{\op{KAC}}(p_\be(U_1),p_\be(U_2)).$$
The complement in $m_\infty$ of the range of $\gamma_*$ has exactly two components
$C_1$ and $C_2$, such that $\bp(U_1)
\in C_1$ and $\bp(U_2)\in C_2$, and the pair $(\lambda(C_1),\lambda(C_2))$
is distributed according to the beta distribution with parameters $(\frac{1}{4},\frac{1}{4})$:
$$\E{f(\lambda(C_1),\lambda(C_2))}= \frac{\Gamma(1/2)}{\Gamma(1/4)^2} \int_{0}^1 dt \big(t(1-t)\big)^{-3/4} f(t,1-t),$$
for any non-negative Borel function $f$ on $\R_+^2$. 
\end{theorem}

\proof
We first explain how the loop $\gamma_*$ is constructed. As in the previous
section, write $(g_{U_1,U_2}(r))_{0\leq r\leq \op{d}_\be(U_1,U_2)}$ for the geodesic path
from $p_\be(U_1)$ to $p_\be(U_2)$ in the tree $\T_\be$, whose range is
the segment $\llbracket p_\be(U_1),p_\be(U_2)\rrbracket$. We already noticed that,
conditionally on the triplet $(\be,U_1,U_2)$ the process
$$ \Big(Z_{g_{U_1,U_2}(r)}-Z_{p_{\be}(U_1)}\Big)_{0\leq r \leq d_{\be}(U_1,U_2)},$$ 
is a standard linear Brownian motion. Hence this process a.s. attains its minimal 
value at a unique time $r_0\in]0, d_{\be}(U_1,U_2)[$, and we put $\beta = g_{U_1,U_2}(r_0)$.
Since there are only countably many values of $r\in ]0, d_{\be}(U_1,U_2)[$ such that
$g_{U_1,U_2}(r)$ has multiplicity $3$ in $\T_\be$, it is also clear that 
$\beta$ has multiplicity $2$ in $\T_\be$, a.s. Write $\mathcal{C}^\circ_1$ and $\mathcal{C}^\circ_2$
for the two connected components of $\T_\be\backslash \{\beta\}$, ordered in
such a way that $p_\be(U_1)\in\mathcal{C}^\circ_1$ and $p_\be(U_2)\in\mathcal{C}^\circ_2$, and set
$\mathcal{C}_1=\mathcal{C}^\circ_1\cup\{\beta\}$, $\mathcal{C}_2=\mathcal{C}^\circ_2\cup\{\beta\}$.
Then $\Pi(\mathcal{C}_1)$ and $\Pi(\mathcal{C}_2)$ are closed subsets of $m_\infty$
whose union is $m_\infty$. Furthermore,
the discussion at the beginning of Section 3 of \cite{Geo} shows that the boundary of 
$\Pi(\mathcal{C}_1)$, or equivalently the boundary of $\Pi(\mathcal{C}_2)$, coincides
with the set $\Pi(\mathcal{C}_1)\cap \Pi(\mathcal{C}_2)$ of all points $x\in m_\infty$
that can be written as $x=\Pi(a_1)=\Pi(a_2)$ for some $a_1\in\mathcal{C}_1$ and 
$a_2\in \mathcal{C}_2$. In particular, the interiors of $\Pi(\mathcal{C}_1)$ and of $\Pi(\mathcal{C}_2)$ are disjoint.
Notice that
$\bp(U_1)$ belongs to the interior of $\Pi(\mathcal{C}_1)$, and $\bp(U_2)$ belongs to the interior of $\Pi(\mathcal{C}_2)$, almost
surely: To see this, observe that for almost every (in the sense of the volume measure $\op{Vol}$)
point $a$ of $\T_\be$, the equivalence class of $a$ for $\approx$ is a singleton, and
thus $\Pi^{-1}(\bp(U_1))$ and $\Pi^{-1}(\bp(U_2))$ must be singletons almost surely. 

Since $\beta$ has multiplicity $2$ in $\T_\be$, Theorem 7.6 in \cite{Geo} implies that 
there are exactly two distinct geodesic paths from $\rho_*$ to $\Pi(\beta)$, and that these paths
are simple geodesics in the sense of \cite[Section 4]{Geo}. We denote these geodesic paths by $\phi_1$
and $\phi_2$. From the definition of simple geodesics, one easily gets that
$\phi_1(s)=\phi_2(s)$ for every $0\leq s\leq s_0$, where
$$s_0:= \max\Big( \min_{a\in \mathcal{C}_1} Z_a, \min_{a\in \mathcal{C}_2} Z_a\Big) - \underline Z.$$
Note that $\{\phi_1(s):0\leq s<s_0\}$ is contained in the interior of $\Pi(\mathcal{C}_i)$, where $i\in\{1,2\}$
is determined by the condition $a_*\in\mathcal{C}_i$.
Furthermore, the definition of simple geodesics shows that
$$\Pi(\mathcal{C}_1)\cap \Pi(\mathcal{C}_2)=\{\phi_1(s):s_0\leq s\leq D(\rho_*,\Pi(\beta))\} \cup \{\phi_2(s):s_0\leq s\leq D(\rho_*,\Pi(\beta))\}.$$ 

We define $\gamma_*$ by setting
$$\gamma_*(t)=\left\{ 
\begin{array}{ll}
\phi_1(t)&\hbox{if }0\leq t\leq D(\rho_*,\Pi(\beta)),\\
\phi_2(2D(\rho_*,\Pi(\beta))-t)\qquad &\hbox{if }D(\rho,\Pi(\beta))\leq t\leq 2D(\rho_*,\Pi(\beta)).
\end{array}
\right.$$
Then $\gamma_*$ is a loop starting and ending at $\rho_*$. Furthermore $\gamma_*$ separates $\bp(U_1)$ 
and $\bp(U_2)$, since any continuous path in $m_\infty$ starting from $\bp(U_1)$ will
have to hit the boundary of $\Pi(\mathcal{C}_1)$ before reaching $\bp(U_2)$. Finally the length of 
$\gamma_*$ is
$$L(\gamma_*)= 2 D(\rho_*,\Pi(\beta)) = 2(Z_\beta-\underline Z)=
D(\rho_*,\bp(U_1)) + D(\rho_*,\bp(U_2)) - 2\,\op{d}_{\op{KAC}}(p_\be(U_1),p_\be(U_2)).$$

We next verify that $\gamma_*$ is the unique loop in $S(\bp(U_1),\bp(U_2),\rho_*)$ with minimal length.
Let $\gamma$ be a path  in $S(\bp(U_1),\bp(U_2),\rho_*)$ indexed by the interval $[0,T]$.
The image
under $\Pi$ of the path $g_{U_1,U_2}$ is a continuous path from $\bp(U_1)$ to $\bp(U_2)$, which must 
intersect the range of $\gamma$. Hence the range of $\gamma$ contains at least one
point $y$ such that $y=\Pi(a)$ for some $a\in \llbracket p_\be(U_1),p_\be(U_2)\rrbracket$.
Since $\gamma(0)=\gamma(T)=\rho_*$, we have
$$L(\gamma)\geq 2\,D(\rho_*,y) = 2(Z_a-\underline Z)$$
using property 1 of the distance $D$ in Section \ref{Brcactus}. Since $Z_a\geq Z_\beta$, we thus obtain 
that $L(\gamma)\geq L(\gamma_*)$. 

Let $\tau\in[0,T]$ be such that $y=\gamma(\tau)$. The preceding considerations show that the equality $L(\gamma)=L(\gamma_*)$
can hold only if $a=\beta$ and if furthermore the paths $(\gamma(\tau - t),0\leq t\leq \tau)$ and
$(\gamma(\tau+t),0\leq t\leq T-\tau)$ have length $D(\rho_*,\Pi(\beta))$, so that these paths must
coincide (up to reparametrization) with geodesics from $\Pi(\beta)$ to $\rho_*$. We conclude that any
minimizing path $\gamma$ coincides with $\gamma_*$, up to reparametrization and time-reversal. 

In order to complete the proof of the theorem, we first need to identify the connected components
of the complement of the range of $\gamma_*$ in $m_\infty$. Consider the case when
$a_*$ belongs to $\mathcal{C}_1$, and set
$$\mathcal{R}:= \{\phi_1(s):0\leq s< s_0\} \subset \Pi(\mathcal{C}_1).$$
Write $\op{Int}(\Pi(\mathcal{C}_i))$ for the interior of $\Pi(\mathcal{C}_i)$, for $i=1,2$. Then 
the connected components
of the complement of the range of $\gamma_*$ in $m_\infty$ are
$$C_1=\op{Int}(\Pi(\mathcal{C}_1)) \backslash \mathcal{R}\ ,\ C_2=\op{Int}(\Pi(\mathcal{C}_2)).$$
This easily follows from the preceding considerations: Note for instance that $\op{Int}(\Pi(\mathcal{C}_2))$
is the image under $\Pi$ of a connected subset of $\mathcal{C}_2$, and is therefore connected.
From this identification, we get
\begin{equation}
\label{connectedcompo}
\lambda(C_1)=\op{Vol}(\mathcal{C}_1)\ ,\ \lambda(C_2)=\op{Vol}(\mathcal{C}_2)=1-\op{Vol}(\mathcal{C}_1),
\end{equation}
using the fact that the range of $\gamma_*$ has zero $\lambda$-measure (this can be seen from
the uniform estimates on the measure of balls found in Section 6 of \cite{Geo}). Clearly the same
identities (\ref{connectedcompo}) remain valid in the case when 
$a_*$ belongs to $\mathcal{C}_2$.

To complete the proof, we need to compute the distribution 
of $\op{Vol}(\mathcal{C}_1)$. To this end it will be convenient to use the invariance of the law
of $\T_\be$ under uniform re-rooting (see e.g. \cite{LGW}). Let $U$ be a random variable uniformly distributed 
over $[0,1]$, and let $\alpha$ be the (almost surely unique) vertex of $\llbracket p_\be(0),p_\be(U) \rrbracket$
such that $Z_\alpha=\min_{a\in\llbracket p_\be(0),p_\be(U) \rrbracket}Z_a$. Then, if $\mathcal{C}^\circ$ is the connected
component of $\T_\be\backslash\{\alpha\}$ containing $p_\be(U)$, the invariance of the CRT under uniform re-rooting
implies that
$$\op{Vol}(\mathcal{C}_1)
\build{=}_{}^{\rm(d)} \op{Vol}(\mathcal{C}^\circ).$$
Now notice that conditionally on the pair $(\be,U)$, the
random variable $H=\op{d}_\be(p_\be(0),\alpha)$
is distributed according to the arc-sine law on $[0,\be_U]$, with density
$$\frac{1}{\pi\sqrt{s(\be_U-s)}}.$$
Moreover,
$$\op{Vol}(\mathcal{C}^\circ)=  D_\be(U,\be_U-H)-G_\be(U,\be_U-H)$$
where we use the same notation as in the preceding section, for $r\in]0,1[$
and $u\in[0,\be_r]$,
\begin{align}
\label{separatech}
&G_\be(r,u) = \max\{ s \leq r:\be_s= \be_r-u\},\notag\\ 
&D_\be(r,u) = \min\{s \geq r:\be_s=\be_r-u\}.
\end{align}
From the previous remarks, we have, for any non-negative measurable function $g$ on $[0,1]$,
\begin{equation}
\label{separatech0}
\E{g(\op{Vol}(\mathcal{C}_1))}
=\E{g(\op{Vol}(\mathcal{C}^\circ))}=\E{\int_0^1 ds\int_0^{\be_s} \frac{dh}{\pi\sqrt{h(\be_s-h)}}\;g(D_\be(s,h)-G_\be(s,h))}.
\end{equation}
In order to compute the right-hand side, it is convenient to argue first under the It\^o measure
$n(de)$ of positive excursions of linear Brownian motion (see e.g. Chapter XII of \cite{RY}, where 
the notation $n_+(de)$ is used). Let $\sigma(e)$ denote the duration of excursion $e$, and define $D_e(r,u)$ and $G_e(r,u)$,
for $r\in]0,\sigma(e)[$ and $0\leq u\leq e(r)$, in a way analogous to (\ref{separatech}). Also write
$$q_h(t)= \frac{h}{\sqrt{2\pi t^3}}\,\exp -\frac{h^2}{2t}$$
for the density of the hitting time of $h>0$ by a standard linear Brownian motion. Then, an application of Bismut's decomposition of the
It\^o measure, in the form stated in \cite[Lemma 1]{URT}, 
gives for every non-negative measurable function $f$ on $\R_+^2$,
\begin{align}
\label{separatech1}
&\int n(de)\int_0^{\sigma(e)} ds\int_0^{e(s)} \frac{dh}{\pi\sqrt{h(e(s)-h)}}\;f\left(\sigma(e), D_e(s,h)-G_e(s,h)\right)\notag\\
&\qquad=\int_0^\infty du \int_0^u\frac{dh}{\pi\sqrt{h(u-h)}}\int_0^\infty dt\,q_{2h}(t)\int_0^\infty dt'\,q_{2(u-h)}(t')\,f(t+t',t)\notag\\
&\qquad=\frac{1}{\pi} \int_0^\infty \frac{dh}{\sqrt{h}}  \int_0^\infty \frac{dh'}{\sqrt{h'}} 
\int_0^\infty dt\,q_{2h}(t)\int_0^\infty dt'\,q_{2h'}(t')\,f(t+t',t)\notag\\
&\qquad=\frac{1}{\pi} \int_0^\infty dt\int_0^\infty dt' \,f(t+t',t) \Big(\int_0^\infty \frac{dh}{\sqrt{h}} \,q_{2h}(t)\Big)
 \Big(\int_0^\infty \frac{dh'}{\sqrt{h'}} \,q_{2h'}(t')\Big).
 \end{align}
 We easily compute
 $$\int_0^\infty \frac{dh}{\sqrt{h}} \,q_{2h}(t)= 2^{-3/4}(2\pi)^{-1/2}\Gamma(3/4)\,t^{-3/4}.$$
 Hence, using also the identity $\Gamma(1/4)\Gamma(3/4)= \pi\sqrt{2}$,
 we see that the right-hand side of (\ref{separatech1}) is equal to
 $$\frac{2^{-3/2}}{\Gamma(1/4)^2} \int_0^\infty d\ell\int_0^\ell dt \,f(\ell,t)\,(t(\ell-t))^{-3/4}.$$
 We can condition the resulting formula on $\{\sigma=1\}$, using the fact that the density of $\sigma(e)$ under $n(de)$ is equal to $\frac{1}{2}(2\pi \ell^3)^{-1/2}$, and we
 conclude that 
 \begin{align*}
 &\E{\int_0^1 ds\int_0^{\be_s} \frac{dh}{\pi\sqrt{h(\be_s-h)}}\;g(D_\be(s,h)-G_\be(s,h))}\\
 &=n\Big(\int_0^{\sigma(e)} ds\int_0^{e(s)} \frac{dh}{\pi\sqrt{h(e(s)-h)}}\;g\left(D_e(s,h)-G_e(s,h)\right)
 \;\Big|\; \sigma =1\Big)\\
 &= \frac{\sqrt{\pi}}{\Gamma(1/4)^2} \int_0^1 dt\,(t(1-t))^{-3/4}\,g(t).
 \end{align*}
We now see that the last assertion of the theorem follows from (\ref{separatech0}). 
\endproof

\section{Appendix}

This section is devoted to the proof of the fact, mentioned in Section
\ref{convmapcactus}, that if $\mathbf{q}=(q_1,q_2,\ldots)$ is a
sequence with finite support, such that $q_k>0$ for some $k\geq 3$,
then there exists a constant $a>0$ such that
$a\mathbf{q}=(aq_1,aq_2,\ldots)$ is regular critical in the sense of
\cite{MaMi,MieInvar}. We briefly discuss case (A2), which is easier. Following
\cite{MaMi}, we define
$$f_{\mathbf{q}}(x)=\sum_{k\geq 0}x^k\binom{2k+1}{k}q_{2k+2}\, ,\qquad
x\geq 0\, .$$ 
By \cite[Proposition 1]{MaMi}, the Boltzmann measure $W_{\mathbf{q}}$ defined in
Section \ref{convmapcactus} is a finite measure if and only if the
equation 
\begin{equation}\label{eq:2}
f_{\mathbf{q}}(x)=1-\frac{1}{x}\, ,\qquad x>1. 
\end{equation}
has a solution. Since $q_k>0$ for some $k\geq 3$, the function
$f_{\mathbf{q}}$ is a strictly convex polynomial, so there can be
either one or two solutions to this equation. In the first situation,
the graphs of $f_{\mathbf{q}}$ and $x\mapsto 1-1/x$ are tangent at the
unique solution, in which case $\mathbf{q}$ is said to be critical in
the sense of \cite[Definition 1]{MaMi} (it will even be {\em regular} critical in
our case since
$f_\bq(x)$ is finite for every $x>0$). It is then trivial that there
exists a unique $a=a_c>0$ such that the graphs of $f_{a\mathbf{q}}$ and
$x\mapsto 1-1/x$ intersect at a tangency point, and then
$a_c\mathbf{q}$ is regular critical.

Let us turn to case (A1), which is more delicate. 
For every $x,y\geq 0$, we set
\begin{eqnarray*}
f^\bullet_{\mathbf{q}}(x,y)&=&
\sum_{k,k'\geq
 0}x^ky^{k'}\binom{2k+k'+1}{k+1}\binom{k+k'}{k}q_{2+2k+k'}\\
f^{\diamond}_{\mathbf{q}}(x,y)&=&
\sum_{k,k'\geq
 0}x^ky^{k'}\binom{2k+k'}{k}\binom{k+k'}{k}q_{1+2k+k'},
\end{eqnarray*}
defining two convex polynomials in the variables $x$ and $y$.  Proposition 1 of
\cite{MieInvar} asserts that the Boltzmann measure $W_{\mathbf{q}}$ is finite
(then $\mathbf{q}$ is said to be {\em admissible}) if
and only if the equations
\begin{equation}
\label{admissib}
\left\{\begin{array}{ll}
\displaystyle{f^\bullet_{\mathbf{q}}(x,y)=1-\frac{1}{x}}\;,\quad &x>1\\
\noalign{\smallskip}
f^{\diamond}_{\mathbf{q}}(x,y)=y\;,&y>0
\end{array}
\right.
\end{equation}
have a solution $(x,y)$, such that the spectral radius of the matrix
$$M(x,y)=\left(\begin{array}{ccc}
    0 & 0 & x-1\\
    \frac{x}{y}\partial_xf_{\mathbf{q}}^{\diamond}(x,y)
    & \partial_yf^{\diamond}_{\mathbf{q}}(x,y) & 0\\
    \frac{x^2}{x-1}\partial_xf^\bullet_{\mathbf{q}}(x,y) &
    \frac{xy}{x-1}\partial_yf^\bullet_{\mathbf{q}}(x,y) & 0
  \end{array}
\right)$$ is at most $1$. Moreover, a solution $(x,y)$ with
these properties is then unique. 

If the spectral radius of $M(x,y)$ (for this unique solution $(x,y)$) equals $1$,
then we say that $\mathbf{q}$ is critical. It is here even {\em regular}
critical in the terminology of \cite{MieInvar}, since the
functions $f^\bullet_{\mathbf{q}},f^\diamond_{\mathbf{q}}$ are
everywhere finite in our case.
Note that the matrix $M(x,y)$ has nonnegative coefficients, and 
 the Perron-Frobenius theorem ensures that the spectral 
radius of $M(x,y)$ is also the largest real eigenvalue of $M(x,y)$. 
Thus, assuming that ${\bf q}$ is admissible, and letting $(x,y)$
be the unique solution of (\ref{admissib}) such that $M(x,y)$ has
spectral radius bounded by $1$, we see that ${\bf q}$ is
regular critical if and only if $1$ is an eigenvalue of $M(x,y)$,
which holds if and only if the determinant of $\op{Id}-M(x,y)$ vanishes.

For every $x,y>0$, set
$$G(x,y)=f^\bullet_{\mathbf{q}}(x,y)-1+1/x\qquad \mbox{ and
}\qquad H(x,y)=f^\diamond_{\mathbf{q}}(x,y)-y\, .$$ Then $G$ and
$H$ are convex functions on $(0,\infty)^2$.  
A pair $(x,y)\in(0,\infty)^2$ satisfies (\ref{admissib}) if and
only if $G(x,y)=H(x,y)=0$ (notice that the condition $G(x,y)=0$
forces $x>1$). 
The set $\{G=0\}$, resp. $\{H=0\}$ is the
boundary of the closed convex set $C_G=\{G\leq 0\}$, resp. of $C_H=\{H\leq 0\}$, in
$(0,\infty)^2$. 

\begin{lemma}\label{sec:appendix}
{\rm (i)} The set $C_G$ is contained in $(1,\infty)\times(0,A)$, for some $A>0$. 

\noindent{\rm(ii)} The set $C_H$ is bounded.

\noindent{\rm (iii)} If $(x,y)\in C_G$ then
$(x,y')\in C_G$ for every $y'\in(0,y)$.  If $(x,y)\in
C_H$ then $(x',y)\in C_H$ for every $x'\in(0,x)$. 
There exists $\eps>0$ such that $C_H$ does not intersect $[1,\infty)\times(0,\eps)$. 

\noindent{\rm (iv)} For every $a>0$, let $G_a$, resp. $H_a$, be the function
analogous to $G$, resp. to $H$, when $\mathbf{q}$ is replaced 
by $a\mathbf q$. Then $C_{H_a}\subset (0,1]\times (0,\infty)$ for every
large enough $a>0$. Consequently $C_{H_a}\cap C_{G_a}=\varnothing$
for every
large enough $a>0$.
\end{lemma}

\proof (i) This is obvious since $f^\bullet_{\mathbf{q}}(x,y)\geq C\,y^\ell$
for every $x,y>0$, for some constant $C>0$ and some integer $\ell\geq 3$. 

\noindent(ii) Suppose first that there exists an odd integer $\ell\geq 3$
such that $q_\ell >0$. Then, the definition of 
$f^\diamond_{\mathbf{q}}$ shows that there is a positive
constant $c$ such that
$$f^\diamond_{\mathbf{q}}(x,y)\geq c(x^{(\ell-1)/2} + y^{\ell-1}),$$
and it readily follows that $C_H$ is bounded. Consider then the case when there
is an even integer $\ell\geq 4$ such that $q_\ell >0$. Then 
there is a positive constant $c$ such that
$$f^\diamond_{\mathbf{q}}(x,y)\geq c(x^{(\ell-2)/2}y+ y^{\ell-1}),$$
and again this implies that $C_H$ is bounded.

\noindent(iii) The first property is clear since 
$y\mapsto G(x,y)$ is non-decreasing, for every $y>0$.
 Similarly, the second property in (iii) follows from the fact that
$x\mapsto H(x,y)$ is non-decreasing, for every $x>0$. The last property is also clear
since we can find $\eps>0$ such that $f^\diamond_{\mathbf{q}}(x,y)>\eps$
for every $x\geq 1$ and $y>0$ (we use the fact that $\mathbf{q}$ is not supported
on even integers).

\noindent(iv) Suppose first that there there exists an odd integer $\ell\geq 3$
such that $q_\ell >0$. Using the same bound as in the proof of (ii), and
noting that $f^\diamond_{a\mathbf{q}}=a\,f^\diamond_{\mathbf{q}}$, we
see that $H_a(x,y)\leq 0$ can only hold if
$$x^{(\ell-1)/2} + y^{\ell-1} \leq \frac{y}{ca}.$$
It is elementary to check that this implies $x\leq 1$ as soon as
$a$ is large enough. The case when there
is an even integer $\ell\geq 4$ such that $q_\ell >0$
is treated similarly using the bound stated in the proof of (ii). Finally the last
assertion in (iv) follows by using (i). \endproof

Recall that $f^\bullet_{\mathbf{q}}$ and $f^\diamond_{\mathbf{q}}$
are polynomials. It follows that the set $\{G=0\}$ is either empty or a smooth curve
depending on whether the set $\{G\leq 0\}$ is empty or not
(a priori it could happen that  $\{G=0\}=\{G\leq 0\}$ is a singleton, but 
assertion (iii) in the previous lemma shows that this case does not occur).
Similar properties hold for the set $\{H=0\}$.
A simple calculation also shows that
\begin{equation}
\label{determinant}
\op{det}(\op{Id}-M(x,y)) =x^2 \op{det}(\nabla G(x,y), \nabla H(x,y)).
\end{equation}
Consequently, if we assume that $(x,y)$ satisfies (\ref{admissib}),
the condition $\op{det}(\op{Id}-M(x,y))=0$ will hold if and only if the curves 
$C_G$ and $C_H$ are tangent at $(x,y)$.

\begin{figure}[h]
\begin{center}
\includegraphics[width=14cm]{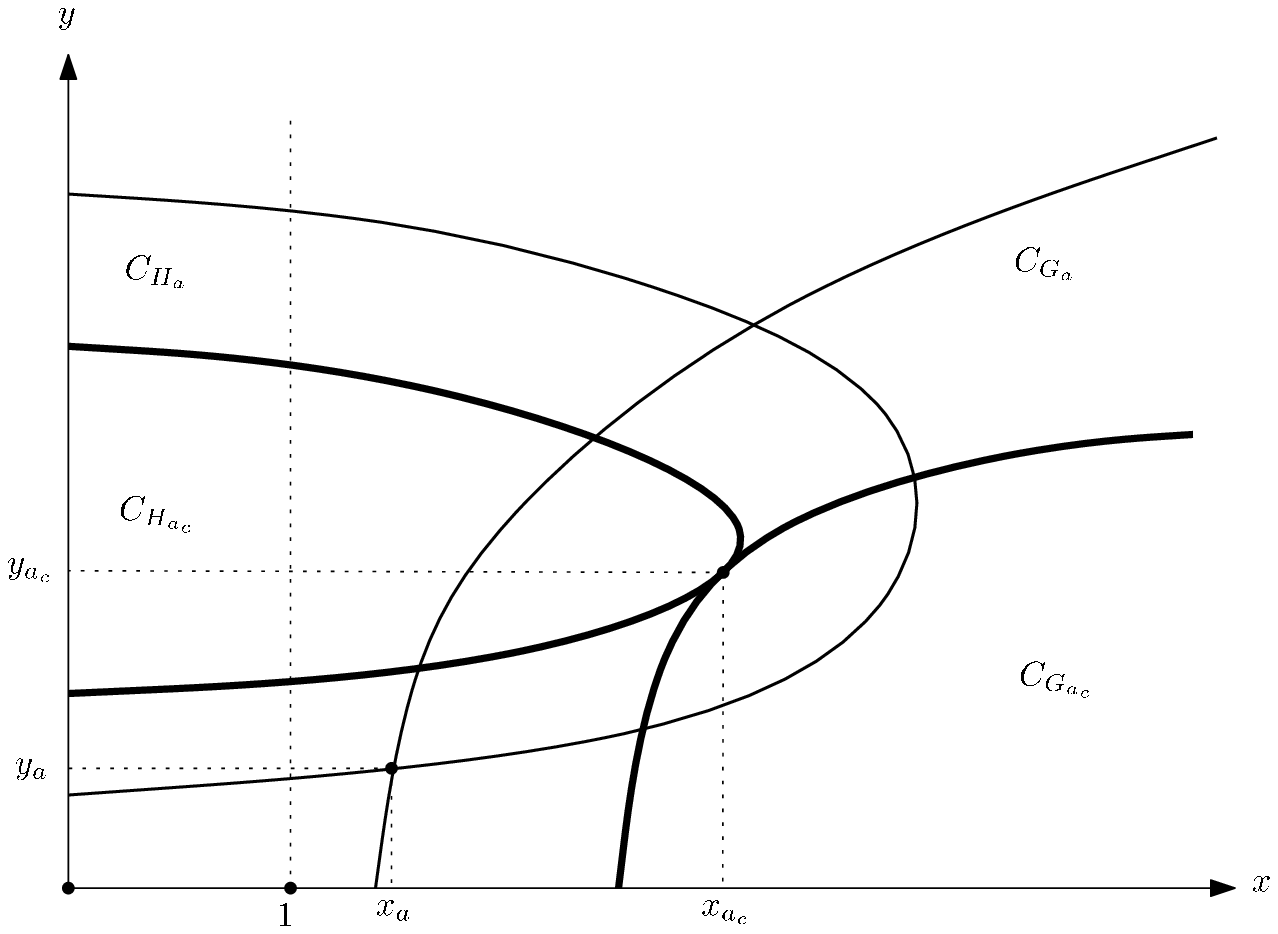}
\caption{Illustration of the sets $C_{G_a}$ and $C_{H_a}$ 
for $0<a<a_c$ and for $a=a_c$}
\end{center}
\end{figure}

\begin{proposition}\label{sec:appendix-1}
Under Assumption (A1),
 there exists a unique positive real $a_c$ such that $a_c\mathbf{q}$
  is regular critical.
\end{proposition}

\proof 
For every $a>0$, write $M_a(x,y)$ for the analogue of the matrix $M(x,y)$
when $\mathbf{q}$ is replaced by $a\mathbf{q}$. Simple counting arguments
(using for instance the BDG bijections and the fact that the sequence
$\mathbf{q}$ has finite support, so that the degrees of faces 
in maps $m$ such that $W_{\mathbf{q}}(m)>0$ are
bounded) show that the Boltzmann measure $W_{a\mathbf{q}}$ is finite for
$a>0$ small enough. Consequently we can fix $a_0>0$ small enough so
that $a_0\mathbf{q}$ is admissible. By previous observations, 
there exists a pair $(x_{a_0},y_{a_0})$ belonging to the intersection
of the curves $\{G_{a_0}=0\}$ and $\{H_{a_0}=0\}$ and such that
the spectral radius of the matrix $M_{a_0}(x_{a_0},y_{a_0})$ is bounded
above by $1$. If the curves $\{G_{a_0}=0\}$ and $\{H_{a_0}=0\}$
are tangent at $(x_{a_0},y_{a_0})$, then (\ref{determinant}) shows that
this spectral radius is equal to $1$, and thus $a_0\mathbf{q}$ is
regular critical.

Suppose that the curves $\{G_{a_0}=0\}$ and $\{H_{a_0}=0\}$
are not tangent at $(x_{a_0},y_{a_0})$.
Then, convexity arguments, using properties (i)--(iii) in Lemma \ref{sec:appendix},
show that the intersection of  $\{G_{a_0}=0\}$ and $\{H_{a_0}=0\}$
consists of exactly two points $(x_{a_0},y_{a_0})$ and $(x'_{a_0},y'_{a_0})$.
By (\ref{determinant}) and the fact that the spectral radius 
of $M_{a_0}(x_{a_0},y_{a_0})$ is bounded
above by $1$, we have 
$$\op{det}(\nabla G_{a_0}(x_{a_0},y_{a_0}), \nabla H_{a_0}(x_{a_0},y_{a_0}))>0,$$
and simple geometric considerations show that $(x_{a_0},y_{a_0})$ must be 
the ``first'' intersection point of $\{G_{a_0}=0\}$ and $\{H_{a_0}=0\}$,
in the sense that $x_{a_0}\leq x'_{a_0}$ and $y_{a_0}\leq y'_{a_0}$. 

Note that both sets $G_a$ and $H_a$ are decreasing functions of $a$, and 
vary continuously with $a$ (as long as they are non-empty).
Geometric arguments, together with property (iv)
of Lemma \ref{sec:appendix}, show that there exists a critical value $a_c>a_0$
such that for $a_0\leq a<a_c$ the curves $\{G_{a}=0\}$ and $\{H_{a}=0\}$
intersect at exactly two points, denoted by $(x_a,y_a)$ and $(x'_a,y'_a)$,
such that $x_a\leq x'_a$ and $y_a\leq y'_a$, and furthermore
the curves $\{G_{a_c}=0\}$ and $\{H_{a_c}=0\}$ are tangent at a
point denoted by $(x_{a_c},y_{a_c})$. Moreover the mapping
$a\mapsto (x_a,y_a)$ is continuous on $[a_0,a_c]$. It follows that the spectral
radius of $M_a(x_a,y_a)$ remains bounded above by $1$ for 
$a\in[a_0,a_c)$: If this were not the case, this spectral radius would take 
the value $1$ at some $a_1\in(a_0,a_c)$ but then by (\ref{determinant}) 
the curves $\{G_{a_1}=0\}$ and $\{H_{a_1}=0\}$ would be tangent
at $(x_{a_1},y_{a_1})$, which is a contradiction. Finally by
letting $a\uparrow a_c$ we get that the spectral radius of 
$M_{a_c}(x_{a_c},y_{a_c})$ is bounded above by $1$, hence
equal to $1$ by (\ref{determinant}) and the fact that $\{G_{a_c}=0\}$ and $\{H_{a_c}=0\}$ are tangent
at $(x_{a_c},y_{a_c})$. We conclude that $a_c\mathbf{q}$ is regular critical. 

The uniqueness of $a_c$ is clear since we can start the previous argument
from an arbitrarily small value of $a_0$ and since the curves 
$\{G_{a}=0\}$ and $\{H_{a}=0\}$ will not intersect when $a>a_c$.
\endproof

 \end{document}